\documentclass[review]{elsarticle}

\usepackage{lineno,hyperref}
\modulolinenumbers[5]

\journal{Journal of \LaTeX\ Templates}

\usepackage{graphicx}
\usepackage{amsmath,amssymb}
\usepackage{xcolor}
\usepackage{mathrsfs}
\usepackage{paralist}
\usepackage{multirow}
\usepackage[T1]{fontenc}
\usepackage[utf8]{inputenc}
\usepackage[english]{babel}
\usepackage{subfigure}
\DeclareMathAlphabet\mathbfcal{OMS}{cmsy}{b}{n}
\usepackage{caption}
\usepackage{lineno}
\usepackage{bbm}
\usepackage[normalem]{ulem}%\usepackage{cancel}
\usepackage[percent]{overpic}

\newcommand{\ud}{\mathrm{d}}

\newcommand{\Reyn}{\mathrm{Re}} %\mbox{\textit{Re}}}
\newcommand{\Ma}{\mathrm{Ma}} %\mbox{\textit{Ma}}}
\newcommand{\pd}[2]{\frac{\partial{#1}}{\partial{#2}}} %partial derivative
\newcommand{\der}[2]{\frac{\ud{#1}}{\ud{#2}}} %derivative
\begin{document}

\begin{frontmatter}

\title{Non-intrusive reduced order models for the accurate prediction of bifurcating phenomena in compressible fluid dynamics}

\author[mymainaddress]{Niccol{\`o} Tonicello\corref{mycorrespondingauthor}}
\cortext[mycorrespondingauthor]{Corresponding author}
\ead{niccolo.tonicello@gmail.com}
\author[mymainaddress]{Andrea Lario}
\author[mymainaddress]{Gianluigi Rozza}
\author[mysecondaryaddress,mythirdaddress]{Gianmarco Mengaldo}

\address[mymainaddress]{Scuola Internazionale Superiore di Studi Avanzati (SISSA), Italy}
\address[mysecondaryaddress]{National University of Singapore, Singapore}
\address[mythirdaddress]{Department of Aeronautics, Imperial College, London, UK}

\begin{abstract}
The present works is focused on studying bifurcating solutions in compressible fluid dynamics. On one side, the physics of the problem is thoroughly investigated using high-fidelity simulations of the compressible Navier-Stokes equations discretised with the Discontinuous Galerkin method. On the other side, from a numerical modelling point of view, two different non-intrusive reduced order modelling techniques are employed to predict the overall behaviour of the bifurcation. Both approaches showed good agreement with full-order simulations even in proximity of the bifurcating points where the solution is particularly non-smooth.
\end{abstract}


\end{frontmatter}
%%%%%%%%%%%%%%%%%%%%%%%%%%%%%%%%%%%%%%%%%%%%%
%%%%%%%%%%%%%%%%%%%%%%%%%%%%%%%%%%%%%%%%%%%%%
%%%%%%%%%%%%%%%%%%%%%%%%%%%%%%%%%%%%%%%%%%%%%
\section{Introduction}
%%%%%%%%%%%%%%%%%%%%%%%%%%%%%%%%%%%%%%%%%%%%%
%%%%%%%%%%%%%%%%%%%%%%%%%%%%%%%%%%%%%%%%%%%%%
%%%%%%%%%%%%%%%%%%%%%%%%%%%%%%%%%%%%%%%%%%%%%
During the last decades, the increasing computational power of modern hardware architectures has made computational fluid dynamics (CFD) a commonly used tool in many different fields of engineering and applied sciences \cite{slotnick2014cfd}. In particular, a large amount of research has been focused on the development of innovative numerical schemes to discretise the partial differential equations (PDEs) which describe the fluid dynamics of many flows of engineering interest. In the context of compressible flows, a large variety of techniques have been developed over the years, from classical Finite Volume (FV), Finite Difference (FD) and Finite Element (FE) methods to more recent high-order schemes such as Discontinuous Galerkin (DG)~\cite{hesthaven2007nodal,cockburn:98,cockburn:98b} or Essentially Non-Oscillatory methods (ENO)~\cite{shu1988efficient,liu1994weighted, shu1998essentially}. All of these schemes are characterised by their own peculiar advantages and drawbacks, making them more or less suitable for the simulation of specific flows. Along these lines, the Discontinuous Galerkin method, which will be extensively employed in this work, was found particularly suitable for the simulation of advected-dominated compressible flows~\cite{Zhang2011PositivitypreservingHO,persson2013shock, mengaldo2015discontinuous,hillewaert2016assessment,moura2017setting, moxey2017towards,moura2017eddy,mengaldo2018spatial_1, mengaldo2018spatial_2,fernandez2018ability,winters2018comparative, fernandez2019non,moura2020viscous,mengaldo2021industry,moura2022spectral}.
Despite the increased efficiency and accuracy of advanced numerical schemes for fluid dynamics, only few applications to cases of practical industrial interest can be found in literature and in most applications the straightforward adoption of such techniques is still unfeasible and cannot adequately answer to the constantly evolving demands of modern industry.
This is the case, for example, when many different parametric configurations need to be numerically simulated during the design process (optimisation), or when rapid input-output responses are required (real-time control). In these situations, classical approaches based on the numerical resolution of partial differential equations modelling the flow are too expensive. On one side, high-fidelity simulations of complex flows are more and more imperative for many specific applications, but, at the same time, the computational cost to obtain accurate solutions is usually too high for most design industrial processes. An expanding branch of research in computational fluid dynamics is currently focused on the use of large databases of high-fidelity simulations to accurately predict a wider range of flow regimes, geometric configurations and parametric dependencies in a more efficient way, by trading computational speed with a little loss in accuracy.
A popular way to reduce the computational cost of high-fidelity simulations is the use of Reduced Order Modelling (ROM) techniques ~\cite{hesthaven2016certified,quarteroni2015reduced,chinesta2016model, benner2017model}. Reduced Order Models are based on the assumption that the dynamical model is primarily driven by a low number of dominant modes. Once these modes are properly identified, a large set of numerical solutions in the parametric space can then be obtained as a combination of the dominant modes.
Of course, such a broad idea can find suitable applications in many different fields, making the use of ROMs a rapidly expanding practice in applied sciences.
The most classical approaches for model reduction are based on standard Galerkin projection on functional spaces. Both the numerical solution and the governing equations themselves are projected on the reduced space, where the reduced space is obtained from Proper Orthogonal Decomposition (POD)~\cite{hotelling1933analysis,lumley1967structure, burkardt2006pod,ballarin2015supremizer,maulik2021pyparsvd,hess2022reduced} of a series of full-order solutions (commonly called \emph{snapshots}). In this case, the knowledge of the full-order model both in terms of numerical scheme (either FE methods or FV methods) and continuum equations (Navier-Stokes equations), is used in the reduced order model. These kind of approaches where a-priori knowledge on the full-model is needed for the ROM are commonly known as \emph{intrusive} ROMs. More recently, spectral POD (or SPOD) has also been used in the context of Reduced-Order modelling~\cite{towne2018spectral,schmidt2019spectral,mengaldo2021pyspod,lario2022neural}.
Looking for more general approaches for complex physics, recent efforts in the ROM community were focused on techniques where the information regarding the numerical scheme or even the equations governing the dynamics of the full-order model is not needed. These approaches are instead commonly referred to as \emph{non-intrusive} ROMs. 
The present work will be focused on the use of non-intrusive ROMs for the accurate prediction of bifurcating solutions of compressible flows, which will be herein modelled using the fully-compressible Navier-Stokes equations.

Reduced Order Models for approximating compressible flows is a research hotspot and both intrusive and non-intrusive approaches have been proposed in literature.
In particular, the latter are preferred in this work due to the tendency of intrusive ROMs to lack robustness when compressibility effects arise. Nevertheless, some works in the literature propose techniques able to mitigate this issue either by employing an empirical quadrature procedure (EQP)~\cite{yano2019} or by rotating the projection subspace, thus introducing the minimal amount of viscosity required to stabilize the system~\cite{BALAJEWICZ2016224}.
On the other hand, non-intrusive ROMs for compressible flows are based on the same principles on which their intrusive counterparts are constructed, i.e. a combination of POD and regression methods \cite{WU2020112766}.

Whereas previous work on ROMs for bifurcating solutions was focused on solid mechanics~\cite{pichi2020reduced,pichi2019reduced}, incompressible fluid dynamics~\cite{pichi2022driving,pitton2017computational,pitton2017application,martin2020reduced,hess2022sparse} and Fluid-Structure Interaction \cite{khamlich2021model}, the use of such techniques for fully-compressible applications still represents a novelty in the field.
In order to deeply investigate the predictive capabilities of ROMs for bifurcation problems of compressible fluid dynamics, sufficiently general, but yet simple test cases are needed. Along these lines, a compressible version of the classical Coanda effect has been considered in the present work~\cite{tritton2012physical,ahmed2019coanda}. The physics of the problem contains rich flow dynamics features with non-negligible compressibility effects, while involving a simple geometry. The problem has been extensively studied in incompressible conditions~\cite{drikakis1997bifurcation,allery2004application,saha2020bifurcation,haffner2020unsteady}, providing a solid validation database in the low-Mach regime. Consequently, the main focus of the present work is to investigate the dependency of the bifurcation phenomenon on the levels of compressibility of the flow, which is quantified through the Mach number. The present work can have a relevant impact in both the physical understanding of compressible bifurcations caused by geometric constraints, and, at the same time, in evaluating the predictive capabilities of non-intrusive ROMs applied to these problems.
The paper will be organised as follows. In section \ref{sec:2} the full-order model is first introduced: both the numerical scheme based on the Discontinuous Galerkin method and its use to discretise the compressible Navier-Stokes equations will be properly discussed. In section \ref{sec:3} the use of non-intrusive Reduced-Order models for the prediction of bifurcating phenomena will be presented. Section \ref{sec:4} will be focused on the numerical simulation of the Coanda effect in sudden expanding bi-dimensional channel flows.  The first part of the section will be focused on the validation of the numerical simulation. The second part, instead, will be dedicated to a physical investigation of the present problem, both in terms of compressibility effects on the bifurcation and its modelling using non-intrusive reduced order models. Conclusions follow in Section \ref{sec:5}.
%%%%%%%%%%%%%%%%%%%%%%%%%%%%%%%%%%%%%%%%%%%%%
%%%%%%%%%%%%%%%%%%%%%%%%%%%%%%%%%%%%%%%%%%%%%
%%%%%%%%%%%%%%%%%%%%%%%%%%%%%%%%%%%%%%%%%%%%%
\section{Full-order model}\label{sec:2}
%%%%%%%%%%%%%%%%%%%%%%%%%%%%%%%%%%%%%%%%%%%%%
%%%%%%%%%%%%%%%%%%%%%%%%%%%%%%%%%%%%%%%%%%%%%
%%%%%%%%%%%%%%%%%%%%%%%%%%%%%%%%%%%%%%%%%%%%%
The fully-compressible Navier-Stokes equations for an ideal gas can be written in a conservative form as:
\begin{equation}
\mathcal{NS}(\textbf{w}) = \frac{\partial \textbf{w}}{\partial t} + \nabla \cdot [ \textbf{F}_{\rm{c}}(\textbf{w}) - \textbf{F}_{\rm{v}}( \textbf{w}, \nabla \textbf{w})] = 0 \quad \quad \forall  \textbf{x}\in \Omega, \;\; t > 0,
\label{EQN_massmomeng}
\end{equation}
with initial and boundary conditions given by:
\begin{align}
%\begin{split}
\textbf{w}(\textbf{x}, 0) & = \textbf{w}_{0}(\textbf{x})  & \forall  &\textbf{x}\in \Omega, \\
\textbf{w}(\textbf{x}, t) & = \textbf{h}_{\rm{g}}               & \forall &\textbf{x}\in \partial \Omega, \;\; t > 0,
%\end{split}
\end{align}
where, $\Omega$ represents the spatial domain and $\partial \Omega$ denotes the domain boundaries. The vector, $\textbf{w} = (\rho, \rho \textbf{u}, \rho E)^{\rm{T}}$ indicates the vector of conserved variables, where $\rho$ is the density, $\textbf{u}$ is the velocity vector and $E$ is the total energy. Furthermore, $\nabla \textbf{w}$ denotes the gradient tensor of the conserved variables while $\textbf{F}_{\rm{c}}$, $\textbf{F}_{\rm{v}} \in \mathbb{R}^{(2+d) \times d} $ are the convective and viscous fluxes respectively (with $d$ denoting the spatial dimensionality of the problem). Their components are given by:
\begin{equation}
\textbf{F}_c(\textbf{w}) =
\left(
\begin{array}{c}
\rho \textbf{u}^{\rm{T}} \\ \rho( \textbf{u} \otimes \textbf{u} ) + p \textbf{I} \\ ( \rho E + p)\textbf{u}^{\rm{T}}
\end{array}
\right),
\;\;\;
\textbf{F}_v(\textbf{w}, \nabla \textbf{w}) =
\left(
\begin{array}{c}
0 \\  \boldsymbol{\tau} \\ \boldsymbol{\tau}\cdot \textbf{u} - \lambda (\nabla T)^{\rm{T}}
\end{array}
\right),
\end{equation}
where $p$ is the pressure, $\lambda$ is the thermal conductivity of the the fluid and $\textbf{I} \in \mathbb{R}^{d \times d}$ is the $d$-dimensional identity matrix. Finally, $\boldsymbol{\tau}$ is the deviatoric part of the stress tensor which can be written as:
\begin{equation}
\boldsymbol{\tau} = 2\mu \Big[ \textbf{S} - \frac{1}{3}(\nabla \cdot \textbf{u})\textbf{I} \Big],
\end{equation}
where $\mu$ is the dynamic viscosity and $\textbf{S}$ is the strain rate tensor given by:
\begin{equation}
\textbf{S} = \frac{1}{2} \Big[ \nabla \textbf{u} + (\nabla \textbf{u})^{\rm{T}} \Big].
\end{equation}
In addition to equations for mass, momentum and energy collectively defined in Eq.~\eqref{EQN_massmomeng}, we require an equation to couple pressure and internal energy and thus close the system of equations. This is given by the ideal gas law written as:
\begin{equation}
p = (\gamma - 1)\Big( \rho E - \frac{1}{2}\rho \textbf{u}\cdot \textbf{u} \Big),
\end{equation}
where $\gamma$ is the ratio of the specific heats. 
%%%%%%%%%%%%%%%%%%%%%%%%%%%%%%%%%%%%%%%%%%%%%
%%%%%%%%%%%%%%%%%%%%%%%%%%%%%%%%%%%%%%%%%%%%%
%%%%%%%%%%%%%%%%%%%%%%%%%%%%%%%%%%%%%%%%%%%%%
\subsection{Discontinuous Galerkin Method}
%%%%%%%%%%%%%%%%%%%%%%%%%%%%%%%%%%%%%%%%%%%%%
%%%%%%%%%%%%%%%%%%%%%%%%%%%%%%%%%%%%%%%%%%%%%
%%%%%%%%%%%%%%%%%%%%%%%%%%%%%%%%%%%%%%%%%%%%%
Equation \ref{EQN_massmomeng} can be written in a even more compact way as:
\begin{equation}
\pd{\textbf{w}}{t} + \nabla \cdot \textbf{F} = 0.
\label{EQN_massmomeng2}
\end{equation}
Before applying the Discontinuous Galerkin approach to the equation \ref{EQN_massmomeng2}, we first introduce the computational domain $\Omega$ as the partition of $N_{e}$ non-overlapping discrete
elements such that $\Omega= \cup_{n=1}^{N_{e}} \Omega_{n}$. Let us also denote the boundary of the $n$-th element as $\partial \Omega_{n}$. 

Secondly, we introduce the space of the test functions in which the numerical solution will be seek into. In particular, a classical choice for Discontinuous Galerkin scheme consists in the functional space:
\begin{equation}
\mathcal{V}^{\delta}:= \{ \phi \in L^{2}(\Omega): \phi \arrowvert_{\Omega_{n}} \in \mathbb{P}_{P}(\Omega_{n}), \forall \Omega_{n} \},
\end{equation}
where $\mathbb{P}_{P}$ is the space of piece-wise continuous polynomial of order not greater than $P$ on $\Omega_{n}$ and 
\begin{equation}
L^{2} (\Omega) = \bigg \{ \phi : \Omega_{n} \rightarrow \mathbb{R} \bigg \arrowvert \int_{\Omega_{n}}|\phi(\textbf{x})|^{2} d \Omega \leq \infty \bigg \}
\end{equation}
with $\mathbb{R}$ being the space of the real numbers. Of course, many different bases can be used to define the functional space $\mathcal{V}^{\delta}$, both modal or nodal. However, to avoid any loss of generality, no specific choice will be considered in the present introduction of the Discontinuous Galerkin method.

The approximation of the global solution $\textbf{w}^{\delta}$ can be defined as:
\begin{equation}
\textbf{w}^{\delta} = \oplus_{n=1}^{N_{e}} \textbf{w}^{\delta}_{n},
\end{equation}
where $\textbf{w}^{\delta}_{n}$ is the local discrete solution: 
\begin{equation}
\textbf{w}^{\delta}_{n} = \sum_{j=0}^{P} \widetilde{\textbf{w}}_{j}(t) \phi_{j}(\textbf{x}).
\end{equation}
The local formulation of the Discontinuous Galerkin method requires $\textbf{w}^{\delta}_{n}$ to satisfy 
\begin{equation}
\int_{\Omega_{n}} \phi_{i} \pd{\textbf{w}^{\delta}_{n}}{t} d\Omega + \int_{\Omega_{n}} \phi_{i} \nabla \cdot \textbf{F}(\textbf{w}^{\delta}_{n}) d\Omega = 0 \quad \forall \phi_{i} \in \mathcal{V}^{\delta}.
\label{weakform}
\end{equation}
The second term in equation \ref{weakform} is the flux term. Upon performing integration by parts, this term can be expressed as
\begin{equation}
\int_{\Omega_{n}} \phi_{i} \nabla \cdot \textbf{F}(\textbf{w}^{\delta}_{n}) d\Omega = - \int_{\Omega_{n}}  \nabla \phi_{i} \cdot \textbf{F}(\textbf{w}^{\delta}_{n}) d\Omega + \oint_{\partial \Omega_{n}} \phi_{i} \widehat{\textbf{F}}(\textbf{w}^{\delta, R}_{n}, \textbf{w}^{\delta, L}_{n}, \widehat{\textbf{n}})d \Gamma,
\end{equation}
where $\widehat{\textbf{n}}$ is the outward-pointing unit normal vector, $(\cdot)^{R}$ and $(\cdot)^{L}$ denote right and left state with respect to the element's interface $\partial \Omega_{n}$ and $\widehat{\textbf{F}}$ is the numerical flux.

After exploiting the form of $\textbf{w}^{\delta}_{n}$, the local discrete weak formulation reads:
\begin{equation}
\sum_{j=0}^{P} \der{ \widetilde{\textbf{w}}_{j}}{\mathrm{t}} \mathrm{M}_{ij} = \mathcal{R}_{i} = \int_{\Omega_{n}}  \nabla \phi_{i} \cdot \textbf{F}(\textbf{w}^{\delta}_{n}) d\Omega - \oint_{\partial \Omega_{n}} \phi_{i} \widehat{\textbf{F}}(\textbf{w}^{\delta, R}_{n}, \textbf{w}^{\delta, L}_{n},\widehat{\textbf{n}})d \Gamma,
\label{weakDG}
\end{equation}
where $M_{ij} = \int_{\Omega_{n}} \phi_{i} \phi_{j} d \Omega$ represents the element local mass matrix and $\mathcal{R}_{i}$ is commonly called \emph{residual}.
Both volume and surface integrals appearing in equation \ref{weakDG} can be either computed analytically or using appropriate quadrature rules. Equation \ref{weakDG} can then be discretised in time either explicitly or implicitly.
If the fluxes depend on the gradients of the conserved variables, like it happens in the Navier-Stokes equations, the problem can be recast as a larger system of equations and a proper definition of numerical fluxes can be formulated.

A large variety of options are available in the open-source code \emph{Nektar}$++$ \cite{cantwell2015nektar++,moxey2020nektar++}, in terms of different set of equations, numerical fluxes, time-integration schemes, $h/p$ adaptation and many other useful features for the simulation of complex fluid dynamics problems \cite{mengaldo2014guide, mengaldo2015dealiasing,mengaldo2015discontinuous,mengaldo2015triple,lombard2016implicit,serson2017direct,nakhchi2020dns,mengaldo2021industry}.
%%%%%%%%%%%%%%%%%%%%%%%%%%%%%%%%%%%%%%%%%%%%%
%%%%%%%%%%%%%%%%%%%%%%%%%%%%%%%%%%%%%%%%%%%%%
%%%%%%%%%%%%%%%%%%%%%%%%%%%%%%%%%%%%%%%%%%%%%
\section{Bifurcation points and reduced order modelling}\label{sec:3}
%%%%%%%%%%%%%%%%%%%%%%%%%%%%%%%%%%%%%%%%%%%%%
%%%%%%%%%%%%%%%%%%%%%%%%%%%%%%%%%%%%%%%%%%%%%
%%%%%%%%%%%%%%%%%%%%%%%%%%%%%%%%%%%%%%%%%%%%%
For the sake of simplicity, any continuous mathematical model we consider in this work can be generally written as a parametric problem of the form:
\begin{equation}
\mathcal{F}(X(\boldsymbol{\mu}),\boldsymbol{\mu}) = 0 \quad \mathrm{ with } \quad \boldsymbol{\mu} \in \mathbb{P}.
\label{eq:par}
\end{equation}
For many applications in applied sciences, computing accurate solutions for many different values of the parameter $\boldsymbol{\mu}$ at low computational expense represents a paramount challenge for industrial applications.
Consequently, from a more mathematical perspective, we have interest in solving the parametric problem \ref{eq:par} for many different values of $\boldsymbol{\mu}$. The ensemble of all the solutions of the reference PDE for each value of the parameter identifies a high-dimensional manifold in the parameter space.
However, for specific problems, multiple solutions could arise from a single specific value of $\boldsymbol{\mu}$. These cases are commonly know as \emph{bifurcations}.

In most of these problems, two separate regions in the parameter space can be identified: one where the solution for any given parameter is uniquely defined and one where multiple solutions (or branches) can be found for any given value of  $\boldsymbol{\mu}$. The specific locations in the parameter space where multiple solutions arise from a single branch are called \emph{bifurcation points}. The dynamics of the system close to these points can be delicate and obtaining accurate approximations of the solution can be a complex task. Consequently, the correct prediction of bifurcating phenomena represents a challenging test for accurate ROM techniques. 

One of the most popular technique is the Reduced Basis method, which allows for a fast evaluation of reliable solutions for new values of the parameter $\boldsymbol{\mu}$. The efficiency of such techniques relies on the decomposition in two phases: the online phase and offline one. The online phase consists in computing full-order solutions for a given set of points in the parameter space. The ensemble of these high-fidelity solutions will be referred to as \emph{snapshots}. The information contained in the pre-computed snapshots can then be used to evaluate new solutions at a significantly lower computational cost. 
The lower costs come from the solution of a smaller problem which essentially consists in a projection into the pre-computed reduced basis spanned by the high-fidelity snapshots. The set of high-fidelity solutions can be expressed as:
\begin{equation}
S = \{ \textbf{X}_{\mathcal{N}}(\boldsymbol{\mu}^{(i)})\} \quad \mathrm{for} \quad i=1,...,N_{s},
\end{equation}
where $N_{s}$ represents the number of snapshots. 

At this stage, it might be reasonable to question how many snapshots are necessary to obtain an accurate and reliable reduced order model. For a giving set of full-order solutions, one way to identify the dominant modes of the problem's dynamics is to use Proper Orthogonal Decomposition (POD). In this way, the optimal reduced basis can be pre-computed and used to evaluate new solutions in the online phase. 
In particular, the POD provides the reduced basis spanned by the first $N$ left singular eigenvectors: $\mathbb{X}_{N} = \mathrm{span} \{\Sigma_{1},...,\Sigma_{N}\}$.
Then, new solutions in the parameter space will be computed as linear combination of the reduced basis as:
\begin{equation}
X_{N} (\boldsymbol{\mu}) = \sum_{n=1}^{N}X_{N}^{(n)}(\boldsymbol{\mu}) \Sigma^{(n)}.
\end{equation}
The final step consists in the computation of the coefficients of the linear combination $X_{N}^{(n)}(\boldsymbol{\mu})$. 
Unfortunately, it is well-known that classical Galerkin projections used during the online phase, for non-linear problems, would still depend on the total number of degrees of freedom of the full-order model. Furthermore, as mentioned in the introduction, classical Galerkin projection for these applications leads to unstable ROMs. Consequently, a series of advanced techniques have been developed over the years to overcome such obstacles. One way, for example, consists in the use of Artificial Neural Networks (ANN) whose main task is to \emph{learn} the projection operation.

In other words, the objective of the ANN is to accurately reconstruct the mapping: 
\begin{equation}
\pi (\boldsymbol{\mu}): \boldsymbol{\mu} \rightarrow X_{N}^{(n)}(\boldsymbol{\mu}).
\end{equation}
The use of ANN can significantly reduce the computational cost of this operation, as well as avoid stability issues.

Such approach has been already successfully applied to ROMs of bifurcating phenomena in fluid dynamics under incompressible conditions. The focus of the present work is to generalise this technique to more complex applications (i.e., for compressible Navier-Stokes equations). The projection step provided by the ANN allows not only to reduce the computational cost but also to decouple the physical origin of the full-order model. The same approach can then be applied to many different variables coming from the full-order model and ideally to many different full-order models in the first place (different set of equations, different schemes). The application of such techniques to compressible fluid dynamics represents the first step towards more extensive multi-physics generalisations of the proposed approach.
In agreement with previous works on reduced order modelling of phenomena dealing with compressible flows~\cite{barone2009stable, mcquarrie2021data}, the ANN procedure is applied to the transformed vector $\textbf{w} = (\frac{1}{\rho},u_{1},u_{2},p)^{T}$ instead of to the conservative variables vector.
%%%%%%%%%%%%%%%%%%%%%%%%%%%%%%%%%%%%%%%%%%%%%
%%%%%%%%%%%%%%%%%%%%%%%%%%%%%%%%%%%%%%%%%%%%%
%%%%%%%%%%%%%%%%%%%%%%%%%%%%%%%%%%%%%%%%%%%%%
\section{Numerical Results}\label{sec:4}
%%%%%%%%%%%%%%%%%%%%%%%%%%%%%%%%%%%%%%%%%%%%%
%%%%%%%%%%%%%%%%%%%%%%%%%%%%%%%%%%%%%%%%%%%%%
%%%%%%%%%%%%%%%%%%%%%%%%%%%%%%%%%%%%%%%%%%%%%

%%%%%%%%%%%%%%%%%%%%%%%%%%%%%%%%%%%%%%%%%%%%%
%%%%%%%%%%%%%%%%%%%%%%%%%%%%%%%%%%%%%%%%%%%%%
%%%%%%%%%%%%%%%%%%%%%%%%%%%%%%%%%%%%%%%%%%%%%
\subsection{Coanda effect}
%%%%%%%%%%%%%%%%%%%%%%%%%%%%%%%%%%%%%%%%%%%%%
%%%%%%%%%%%%%%%%%%%%%%%%%%%%%%%%%%%%%%%%%%%%%
%%%%%%%%%%%%%%%%%%%%%%%%%%%%%%%%%%%%%%%%%%%%%
The Coanda effect can be observed in a wide range of fluid flows: from  aerodynamics problems \cite{wang2009transient,trancossi2011overview,trancossi2011acheon,ahmed2017aerodynamics} to acoustics applications \cite{lubert2010some}, as well as multiphase flows \cite{freire2002bubble} and cardiology studies \cite{ginghina2007coandua}. The present work will be focused on the role played by compressibility in bifurcating phenomena, using a combination of high-fidelity simulations obtained via the Discontinuous Galerkin method and exploring the use of non-intrusive reduced order model to predict the bifurcation. Previous work has been done in the prediction of bifurcations for incompressible flows \cite{hess2019localized,pintore2021efficient,pichi2021artificial} and fluid-structure interaction problems \cite{khamlich2021model}. Ideally, we would like to show that similar non-intrusive algorithms based on Neural Networks can be adopted to the accurate modellisation of bifurcating phenomena in compressible flows.

The present work will be focused on the Coanda effect \cite{durst1974low,fearn1990nonlinear}: a relatively simple test case that can be used to assess the suitability of reduced order modelling techniques for bifurcating solutions in compressible fluid dynamics. From a fluid mechanics standpoint, a recent study by Karantonis et al.\cite{karantonis2021compressibility} was similarly focused on the role played by compressibility for the present configuration, considering, however, only configurations far from the bifurcation point (fully-turbulent case).

The use of a compressible setup for the same problem represents a natural evolution of previous research and, for relatively small Mach numbers, the present implementation can be immediately compared with incompressible flow results for validation purposes. 
%%%%%%%%%%%%%%%%%%%%%%%%%%%%%%%%%%%%%%%%%%%%%
%%%%%%%%%%%%%%%%%%%%%%%%%%%%%%%%%%%%%%%%%%%%%
%%%%%%%%%%%%%%%%%%%%%%%%%%%%%%%%%%%%%%%%%%%%%
\subsubsection{Numerical set-up}
%%%%%%%%%%%%%%%%%%%%%%%%%%%%%%%%%%%%%%%%%%%%%
%%%%%%%%%%%%%%%%%%%%%%%%%%%%%%%%%%%%%%%%%%%%%
%%%%%%%%%%%%%%%%%%%%%%%%%%%%%%%%%%%%%%%%%%%%%
The geometrical set-up of the present test case consists in a simple two-dimensional duct with a sudden area expansion.

In terms of boundary conditions, at the inlet, a parabolic velocity profile and unitary density are imposed, whereas the pressure is left free. On the other hand, at the outflow boundary, the velocity field and density are left free, whereas an outlet pressure is imposed. The choice of the velocity profile and outlet pressure will define the characteristic Mach number, which quantifies the importance of compressibility effects for the present problem. The remaining boundaries are modelled as adiabatic walls. 

\begin{figure}
\centering
\includegraphics[ width=0.9\textwidth]{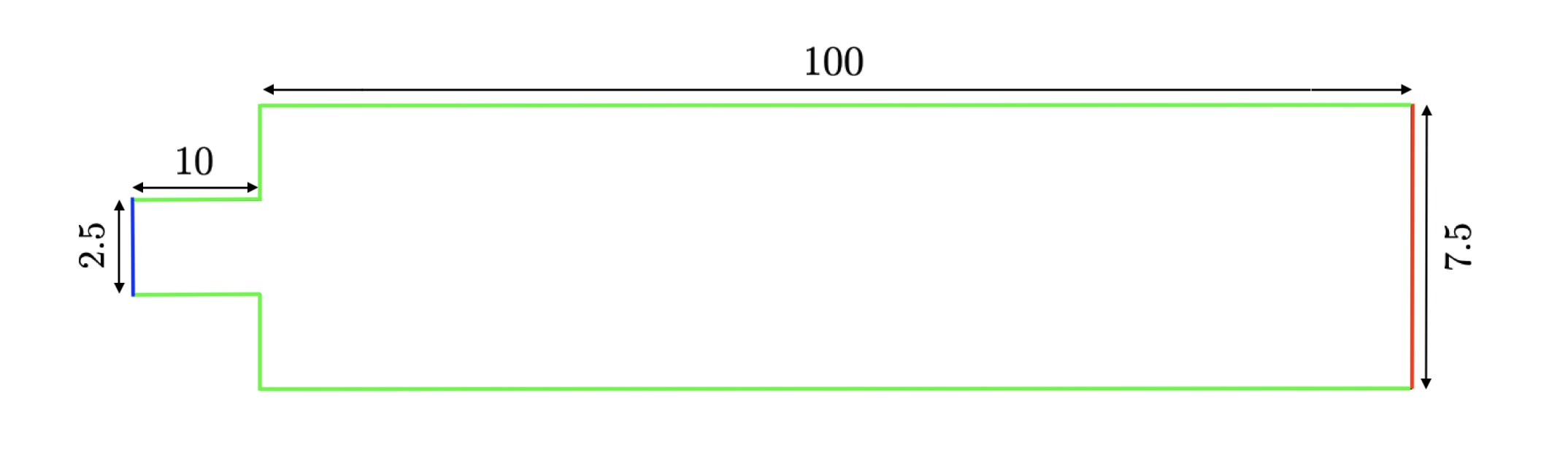}
\caption{Problem's geometry and boundary conditions. Blue, inlet boundary; red, outlet boundary; green, adiabatic wall.}
\label{fig:domain}
\end{figure}

A schematic overview of the problem's geometry and boundary conditions is shown in Figure \ref{fig:domain}.
The present configuration mainly depends on two adimensional numbers: the Reynolds number and the Mach number. Namely, the two quantities are defined as
\begin{equation}
\Reyn = \frac{\rho_{\mathrm{in}} u_{\mathrm{in}} L}{\mu} \quad \mathrm{and} \quad \Ma = \frac{u_{\mathrm{in}}}{\sqrt{\gamma p_{\mathrm{out}} / \rho_{\mathrm{in}}}},
\end{equation}
where $\mu$ is the dynamic viscosity, $\gamma = c_p/c_v$ is the specific heat ratio,  $L$ is the width of the channel at the inlet, and the subscripts $(\cdot)_{\mathrm{in}}$ and $(\cdot)_{\mathrm{out}}$ respectively represent values taken at inlet and outlet boundaries.

For the specific problem herein considered, the same values used by Pichi et al.\cite{pichi2021artificial} have been used. Namely, a parabolic stream-wise velocity profile with maximum peak value of $u_{\mathrm{in}} =20$. For the same reasons, the width of the duct at the inlet has a value of $L=2.5$. Finally, a unitary value has been assigned to density at the inlet. In this way, Reynolds and Mach numbers are entirely defined by the choice of viscosity and outlet pressure respectively.

The first part of the study will be restricted to the low-Mach region in order to provide a suitable framework for comparison with the work by Pichi and co-authors. In this region, the range of viscosity $\mu \in [0.5,2]$ has been considered. In terms of dimensionless numbers, this is equivalent to consider the range of Reynolds numbers $\Reyn \in [25,100]$. A pitchfork bifurcation is expected to take place at $\Reyn \approx 52$ (or $\mu \approx 0.96$). For larger Mach numbers, the location of the bifurcation point is expected to move. This specific matter will be discussed in the second part of this section.

\begin{figure}
\centering
\begin{overpic}[trim=0 210 0 100,clip,width=0.95\columnwidth]{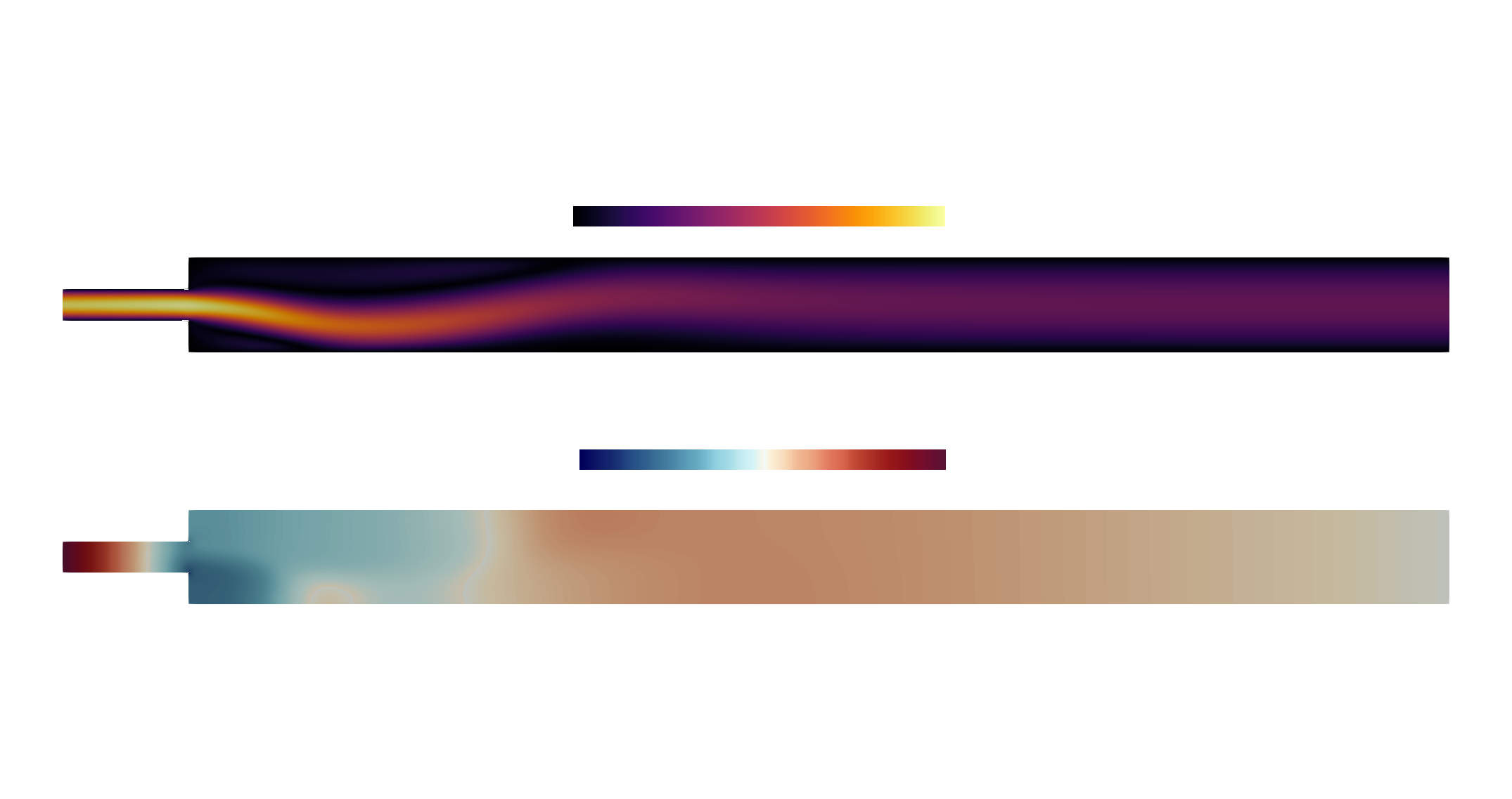}
\put (36.5,30) {$0.0$}
\put (48,32) {$\mathrm{Ma}$}
\put (60.7,30) {$0.48$}
\put (36.5,14) {$0.95$}
\put (48,15) {$p/p_{\mathrm{out}}$}
\put (60.7,14) {$1.11$}
\end{overpic}
\caption{Representative example of asymmetric solution for the case $\mathrm{Ma}=0.3$. Top, Mach number; bottom, pressure field.}
\label{fig:M03ex}
\end{figure}
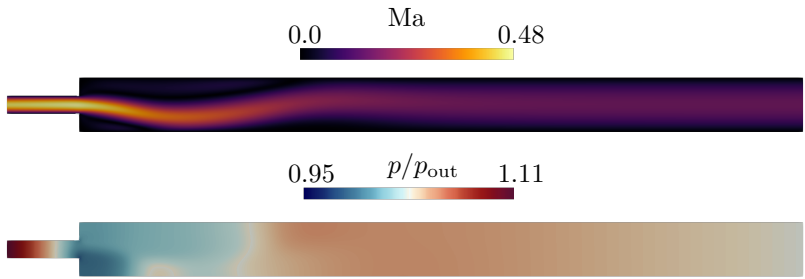

An example of asymmetric bifurcating solution for a relatively small Mach number is shown in Figure \ref{fig:M03ex}. For sufficiently large values of the Reynolds number, the flow field tends to be \emph{attracted} to the solid walls, leading to a symmetry-breaking solution. In the pressure field, it can be noticed that a strong depression region naturally develops on the bottom half of the domain. 
%
%%%%%%%%%%%%%%%%%%%%%%%%%%%%%%%%%%%%%%%%%%%%%
%%%%%%%%%%%%%%%%%%%%%%%%%%%%%%%%%%%%%%%%%%%%%
%%%%%%%%%%%%%%%%%%%%%%%%%%%%%%%%%%%%%%%%%%%%%
\subsubsection{Grid convergence study}
%%%%%%%%%%%%%%%%%%%%%%%%%%%%%%%%%%%%%%%%%%%%%
%%%%%%%%%%%%%%%%%%%%%%%%%%%%%%%%%%%%%%%%%%%%%
%%%%%%%%%%%%%%%%%%%%%%%%%%%%%%%%%%%%%%%%%%%%%
In order to generate high-fidelity snapshots, a grid convergence study has been carried out. Two different flow conditions have been considered. In particular, the outlet pressure has been chosen in order to obtain characteristic Mach numbers equal to $0.3$ and $0.6$. In both cases the value of viscosity was set to $\mu=0.5$. 

A series of subsequently more refined grids has been considered. In particular, the only relevant parameter for the mesh generation can be identified with the number of elements along the wall-normal direction at the inlet. The remaining number of subdivisions of the domain's edges have been chosen in order to obtain identical square elements everywhere in the domain. The number of elements along the wall-normal direction at the inlet varies between $N_{y}=4,6,8,10,12$. 

The numerical simulations herein presented have been obtained using the \emph{Nektar++}\cite{cantwell2015nektar++} compressible flow solver. In particular, a $3^{\mathrm{rd}}$-order Discontinuous Galerkin scheme has been used, coupled with a $3^{\mathrm{rd}}$-order Diagonally Implicit Runge-Kutta \cite{kennedy2016diagonally} time integration . The original Roe solver \cite{toro2013riemann} has been used to compute inviscid numerical fluxes and an interior penalty approach \cite{hartmann2006symmetric} for viscous numerical fluxes.

Time integration is performed till convergence to a stable steady-state solution is reached.

\begin{figure}
\centering
\includegraphics[ width=0.9\textwidth]{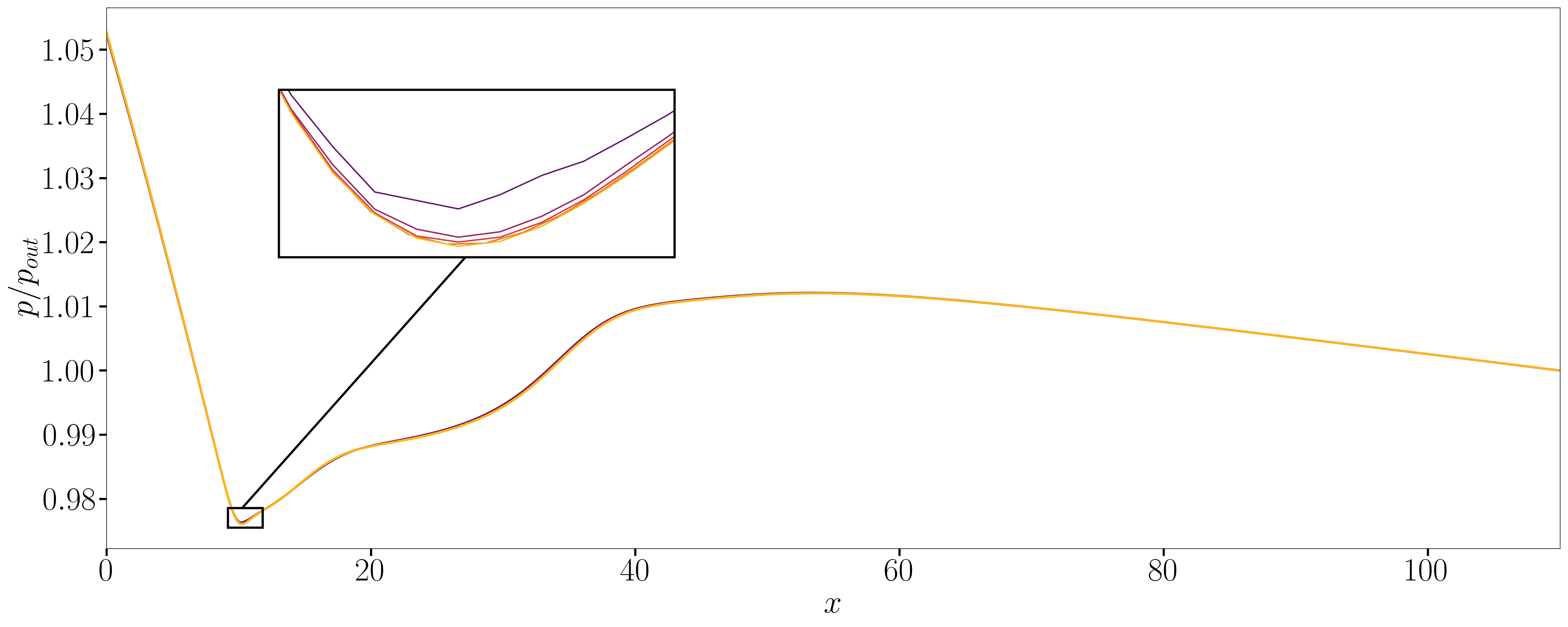}
\caption{Pressure field along the stream-wise direction at $y=0$ for the case $\mathrm{Ma}=0.3$. Color gradient indicates mesh refinement, from the coarsest (purple) to the finest (dark yellow) resolution.}
\label{fig:Ma03x1}
\end{figure}

\begin{figure}
\centering
\includegraphics[ width=0.9\textwidth]{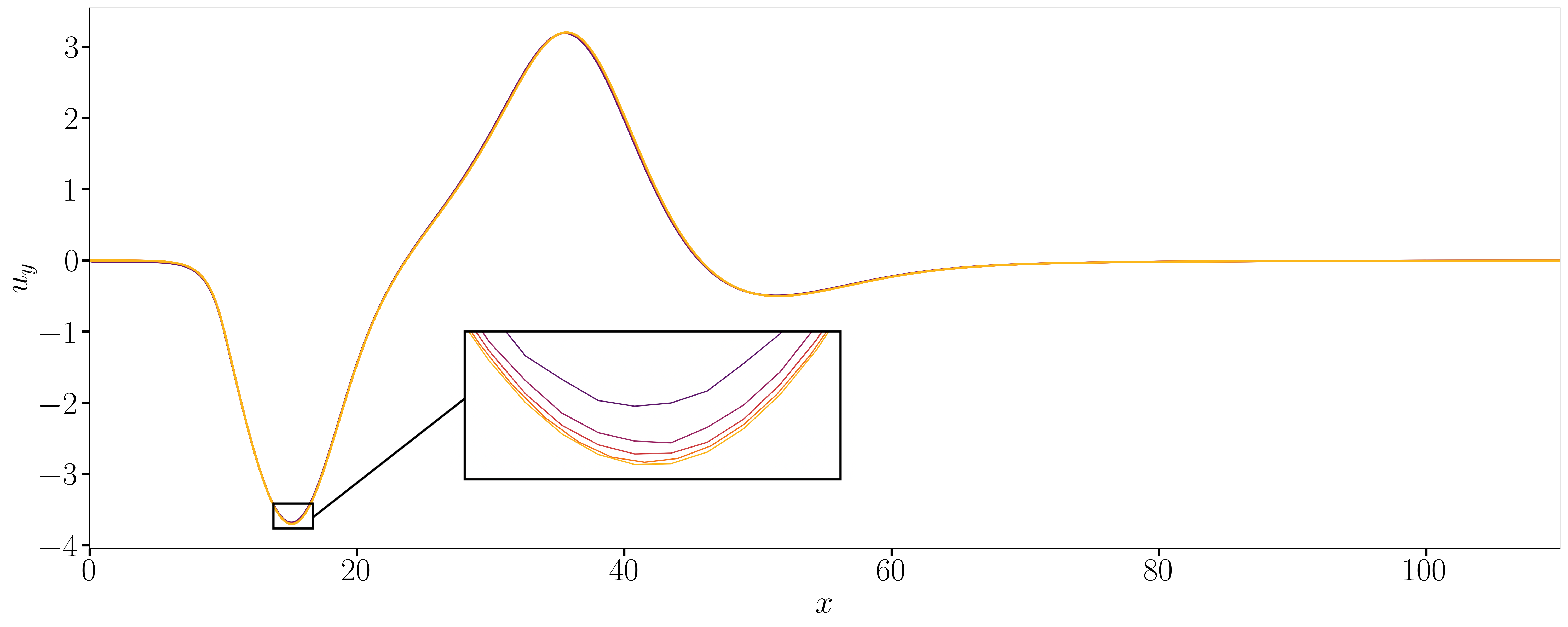}
\caption{Wall normal velocity field along the streamwise direction at $y=0$ for the case $\mathrm{Ma}=0.3$. Color gradient indicates mesh refinement, from the coarsest (purple) to the finest (dark yellow) resolution.}
\label{fig:Ma03x2}
\end{figure}

For the case $\mathrm{Ma}=0.3$, in Figures \ref{fig:Ma03x1} and \ref{fig:Ma03x2}, pressure field and wall-normal velocity profile are respectively shown along the stream-wise direction at the midsection of the channel for different choices of grid size. Both pressure field and wall-normal velocity field are almost identical for different grid sizes. In particular, zooming into the most irregular regions of the flow field, it can be noticed that for resolutions higher than $N_{y}=8$ all the profiles tend to collapse on one another.

\begin{figure}
\centering
\includegraphics[ width=0.9\textwidth]{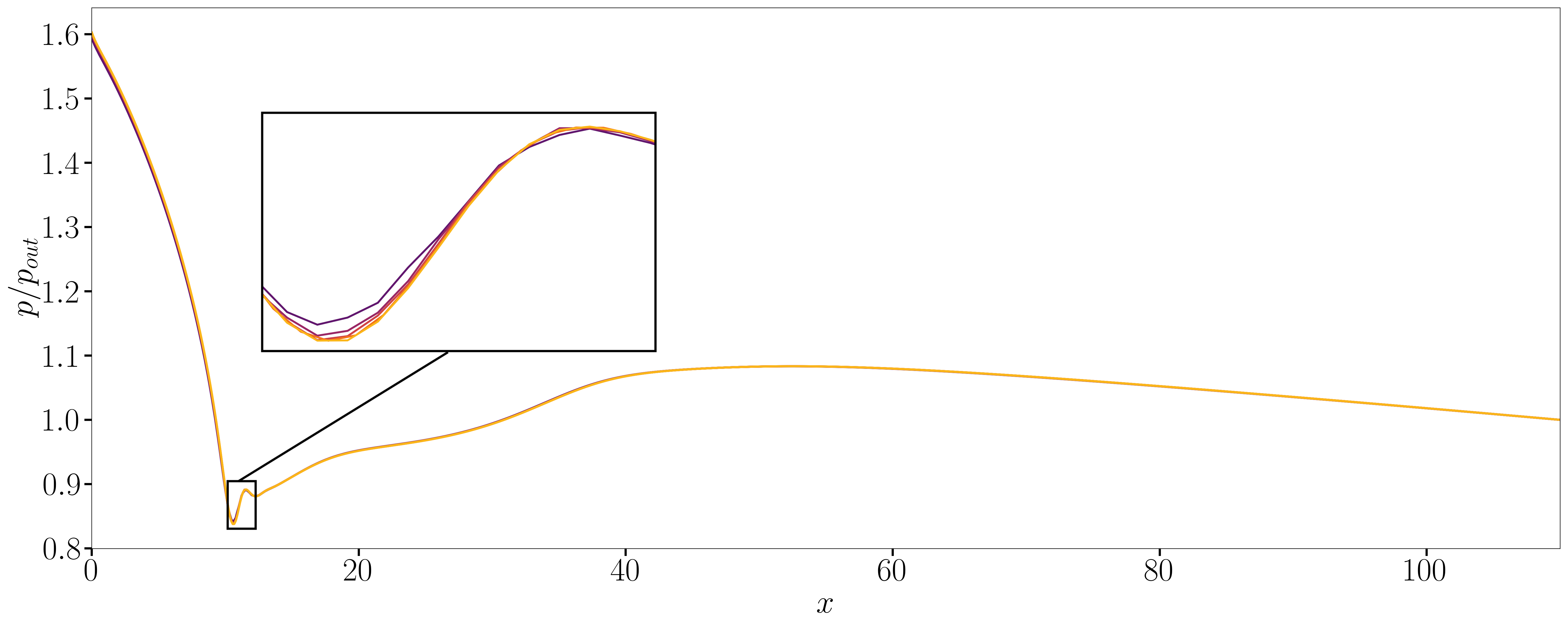}
\caption{Pressure along the streamwise direction at $y=0$ for the case $\mathrm{Ma}=0.6$. Color gradient indicates mesh refinement, from the coarsest (purple) to the finest (dark yellow) resolution.}
\label{fig:Ma06x1}
\end{figure}

\begin{figure}
\centering
\includegraphics[ width=0.9\textwidth]{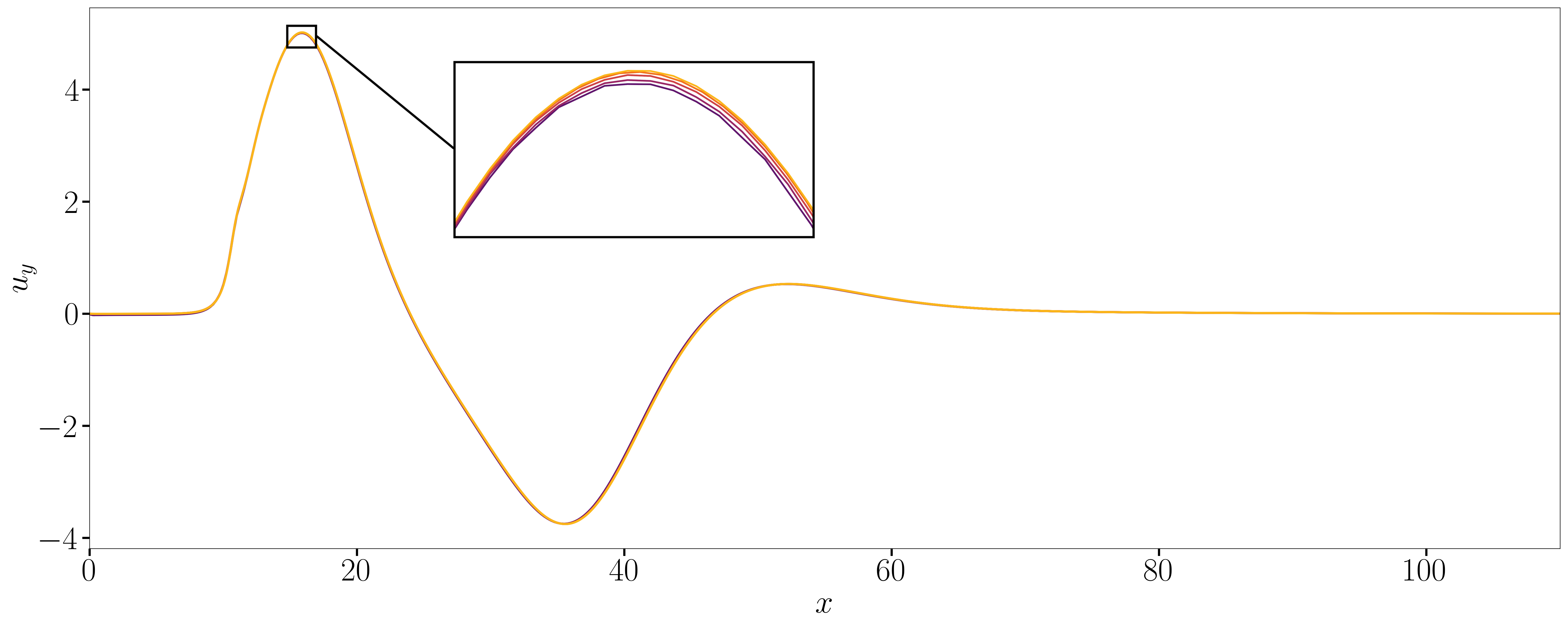}
\caption{Wall normal velocity field along the stream-wise direction at $y=0$ for the case $\mathrm{Ma}=0.6$. Color gradient indicates mesh refinement, from the coarsest (purple) to the finest (dark yellow) resolution.}
\label{fig:Ma06x2}
\end{figure}

The same plots for $\mathrm{Ma}=0.6$ are shown in Figures \ref{fig:Ma06x1} and \ref{fig:Ma06x2} respectively. 


\begin{figure}
\centering
\subfigure[Pressure field.]{\includegraphics[width=0.48\textwidth]{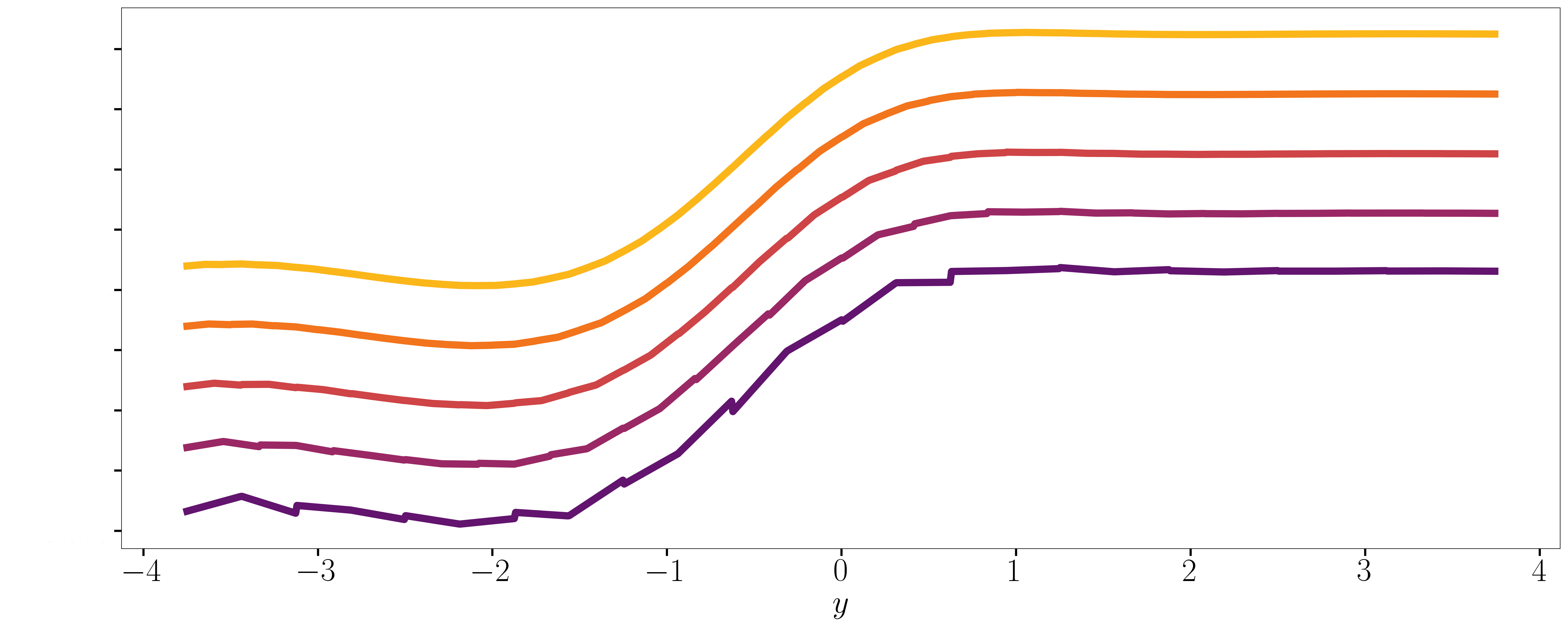}}
\subfigure[Streamwise velocity profile.]{\includegraphics[width=0.48\textwidth]{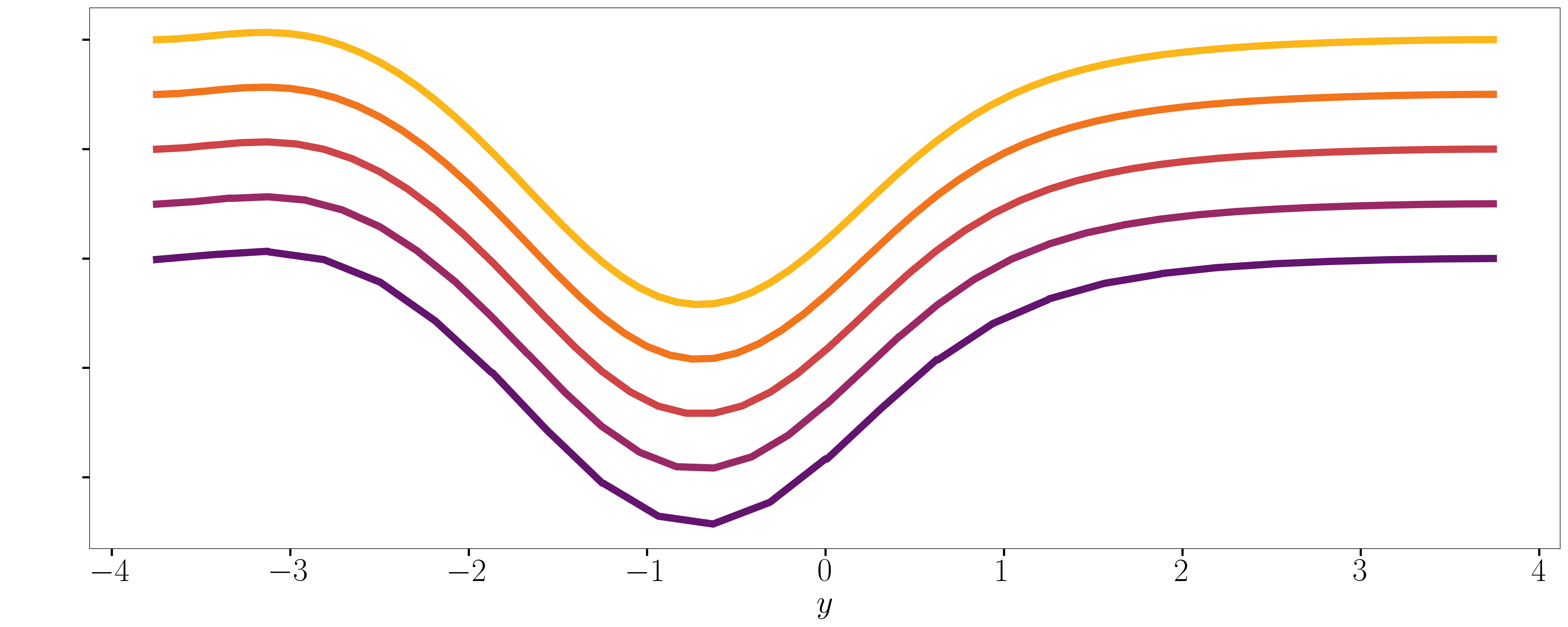}}
\caption{Pressure field (left) and streamwise velocity profile (right) along the wall normal direction at $x=15$ for the case $\mathrm{Ma}=0.3$. Color gradient indicates mesh refinement, from the coarsest (purple) to the finest (dark yellow) resolution. Profiles are shifted for clarity purposes.}
\label{fig:Ma03y15}
\end{figure}

\begin{figure}
\centering
\subfigure[Pressure field.]{\includegraphics[width=0.48\textwidth]{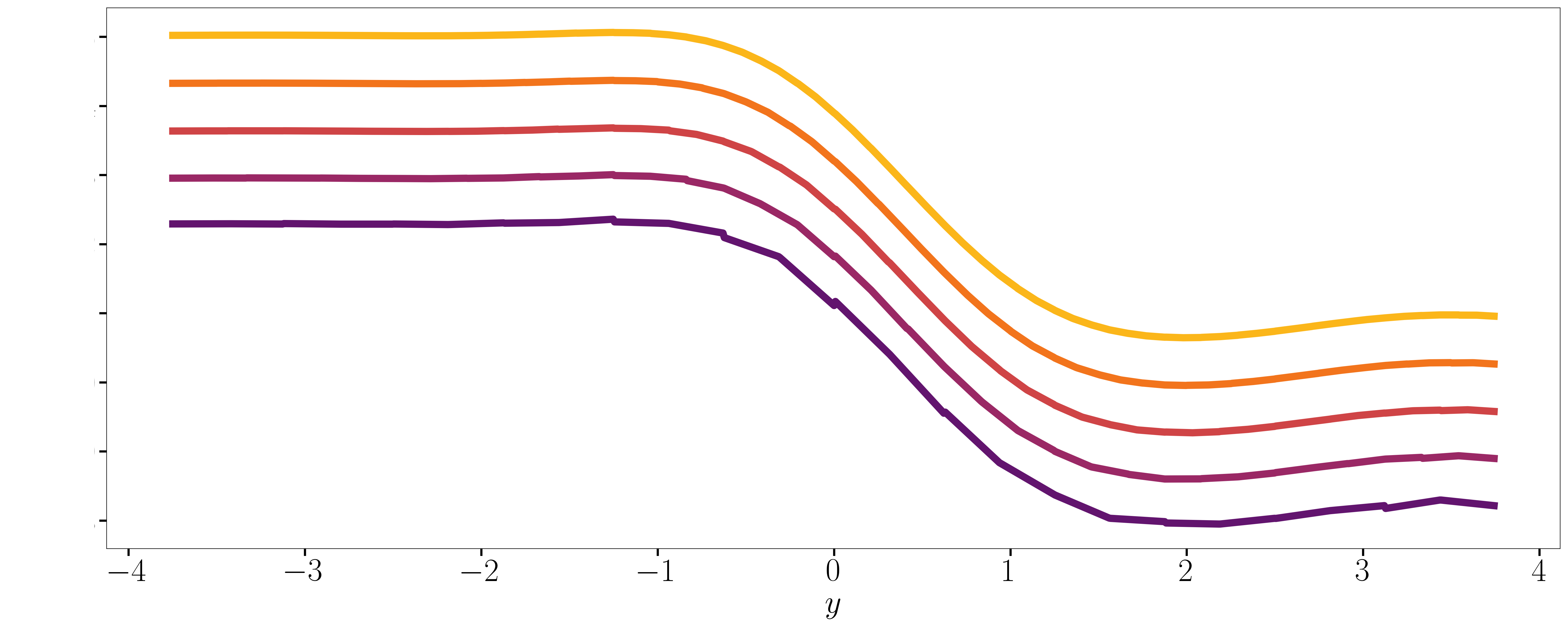}}
\subfigure[Streamwise velocity profile.]{\includegraphics[width=0.48\textwidth]{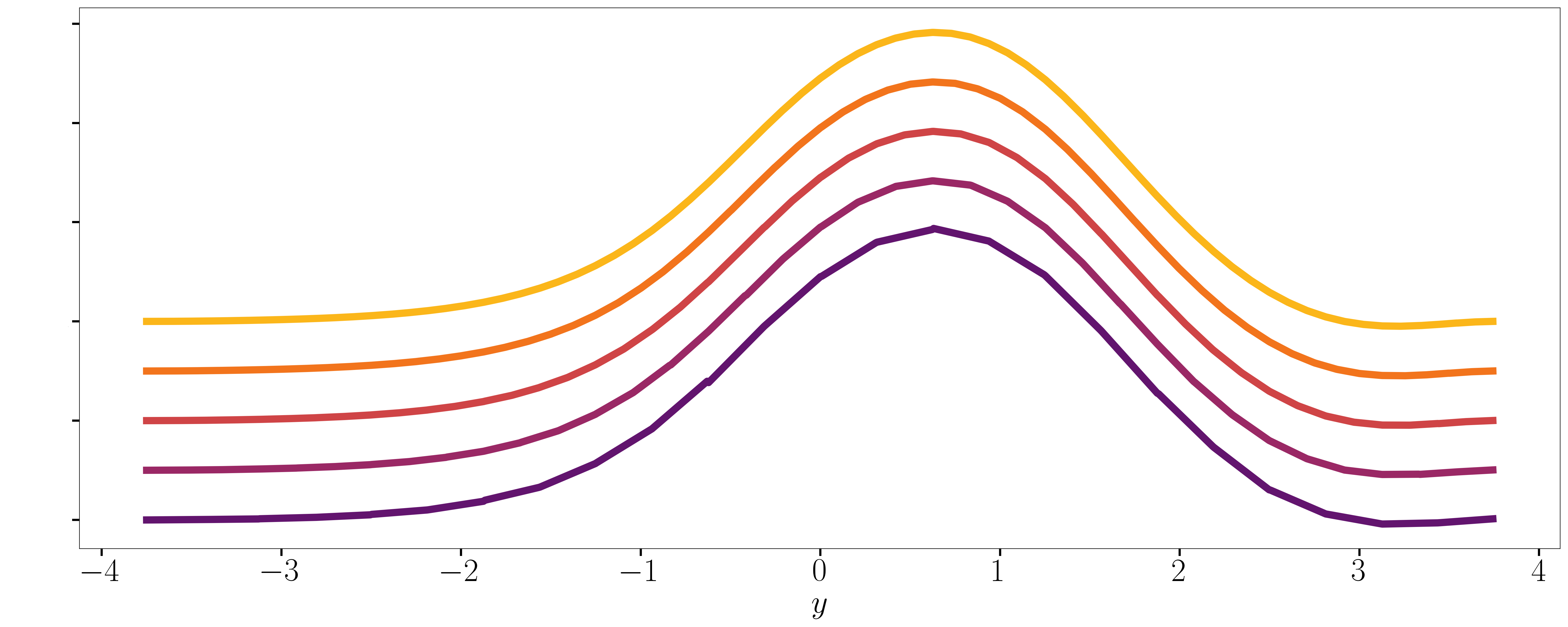}}
\caption{Pressure (left) and streamwise velocity (right) profiles along the wall normal direction at $x=15$ or the case $\mathrm{Ma}=0.6$. Color gradient indicates mesh refinement, from the coarsest (purple) to the finest (dark yellow) resolution. Profiles are shifted for clarity purposes.}
\label{fig:Ma06y15}
\end{figure}

Finally, driven by similar motivations, wall-normal profiles have been considered as well. In particular, pressure field and stream-wise velocity profile have been computed along the wall normal direction at $x=15$ and they are shown in Figures \ref{fig:Ma03y15} and \ref{fig:Ma06y15} for the case $\mathrm{Ma}=0.3$ and $\mathrm{Ma}=0.6$, respectively. 

Similar conclusions between wall-normal and stream-wise profiles can be drawn: for $N_{y}=4$ and $N_{y}=6$ the solution is still relatively irregular, but for smaller grid sizes, the profiles are very similar one another. 
Based on the present preliminary results, we consider the grid characterised by $8$ elements at the inlet to be sufficiently grid-converged to be used as reference for the generation of high-fidelity snapshots. Notice that the total number of degrees of freedom is almost identical to the one used by Pichi et al.\cite{pichi2021artificial}. Namely, in their work, a total of $79868$ degrees of freedom have been used versus the $80352$ degrees of freedom of the present discretisation. Finally, notice that the present discretisation is characterised by a higher order of approximation with respect to the Taylor-Hood finite elements used by Pichi et al.\cite{pichi2021artificial}.
%%%%%%%%%%%%%%%%%%%%%%%%%%%%%%%%%%%%%%%%%%%%%
%%%%%%%%%%%%%%%%%%%%%%%%%%%%%%%%%%%%%%%%%%%%%
%%%%%%%%%%%%%%%%%%%%%%%%%%%%%%%%%%%%%%%%%%%%%
\subsubsection{Full-order model analysis of the bifurcation}
%%%%%%%%%%%%%%%%%%%%%%%%%%%%%%%%%%%%%%%%%%%%%
%%%%%%%%%%%%%%%%%%%%%%%%%%%%%%%%%%%%%%%%%%%%%
%%%%%%%%%%%%%%%%%%%%%%%%%%%%%%%%%%%%%%%%%%%%% 
Before introducing the reduced order techniques used in the present work, a general overview on the dynamics of the bifurcation will be first provided using the full-order model simulations. 
As a first assessment, we want to quantify the relevance of compressibility effects on the overall structure of the bifurcation diagram. 
In particular, five values of the Mach number have been considered: $\mathrm{Ma}=0.3,0.6,0.7,0.8,0.9$.

\begin{figure}
\centering
\begin{overpic}[width=0.9\textwidth]{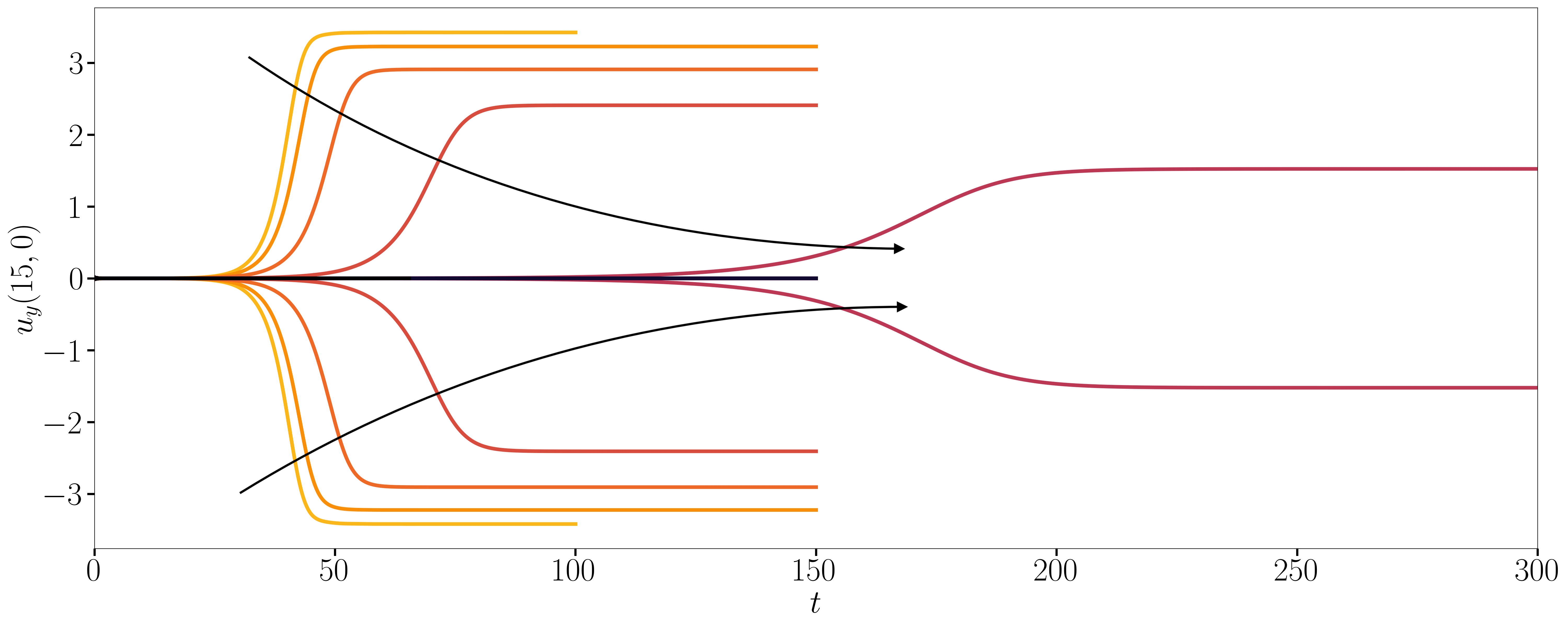}
\put (40,28.5) {$\mu$}
\end{overpic}
\caption{Probed solution at $(x,y)=(15.0,0.0)$ (wall-normal velocity field) for the case $\mathrm{Ma}=0.3$. Color gradient denotes varying values of molecular viscosity, from $\mu=0.5$ (dark yellow) to $\mu=2.0$ (black).}
\label{fig:Ma03timehis}
\end{figure}

In Figure\ref{fig:Ma03timehis}, for the smallest Mach number case, the wall-normal component of the velocity field at  $(x,y)=(15.0,0.0)$ is monitored over time. It can be immediately noticed that for decreasing values of dynamic viscosity (represented by the color gradient) two stable asymmetric branches arise from the unique stable symmetric solution.
It can also be noticed that the simulated time can vary depending on the value of viscosity. For values of dynamic viscosity closer to the bifurcation, in fact, the convergence to one of the two branches is slower as it takes place at longer times.

Using the information contained in the probed solution, it is possible to build the bifurcation diagram using the final value of the wall-normal velocity field at $(x,y)=(15.0,0.0)$.

\begin{figure}
\centering
\includegraphics[width=0.9\textwidth]{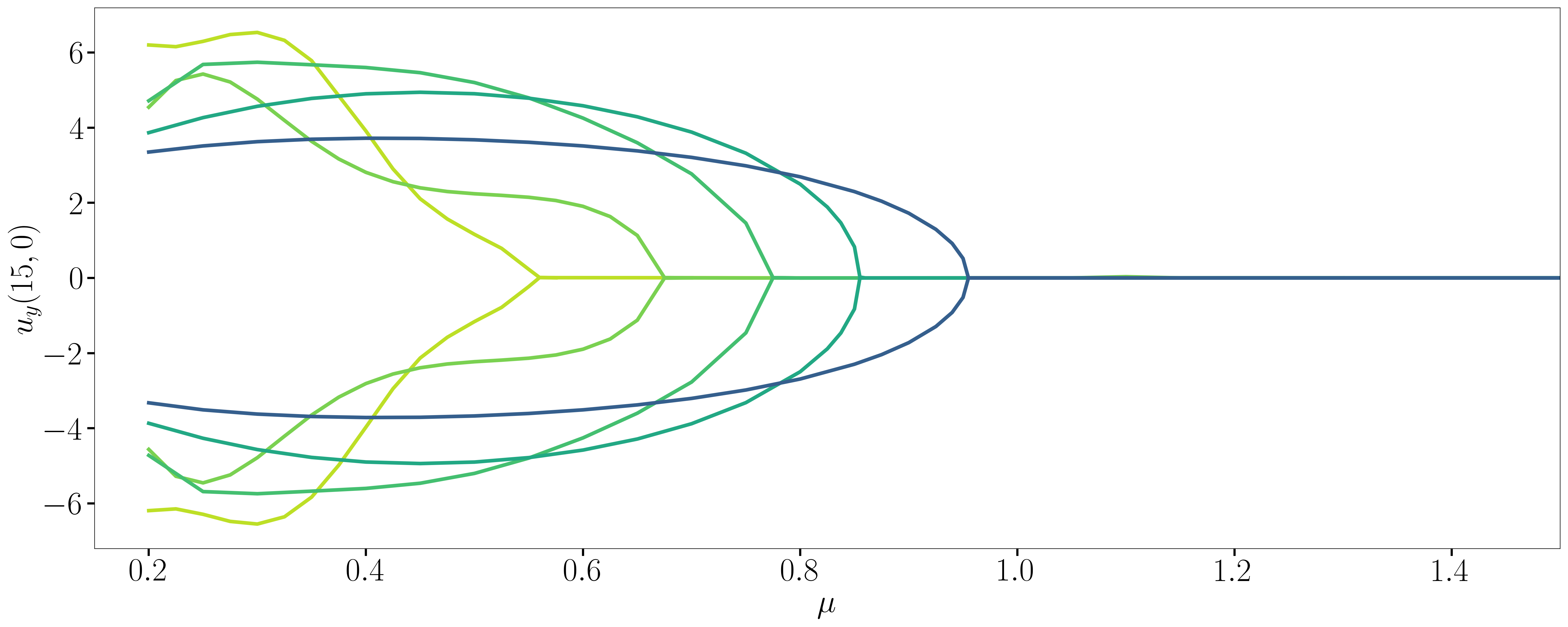}
\caption{Bifurcation diagram using the vertical velocity component probed in the middle of the channel. Color gradients indicates increasing values of the Mach number.}
\label{fig:diagram}
\end{figure}

The approximate bifurcation diagram is reported in Figure \ref{fig:diagram}. It is easy to see that, for small values of the Mach numbers, the bifurcation takes place approximately between $\mu=0.9$ and $\mu=1.0$, as expected for this specific problem. Observing the same diagram for the case $\mathrm{Ma}=0.6$, instead, some differences can spotted. The bifurcation point is located at smaller values of viscosity (between $\mu=0.8$ and $\mu=0.9$). In fact, the numerical solution for $\mu=0.9$ and $\mathrm{Ma}=0.6$ is symmetric, in contrast with the incompressible analogue solution, which is already asymmetric.

For even higher Mach number the overall form of the bifurcation diagram changes significantly and the two asymmetric branches not only arise for smaller valued of viscosity but they are also characterised by a different shape. In proximity of the bifurcation point, the gradients of the bifurcation diagram are smaller, denoting a slower bifurcation with respect to the almost discontinuous behaviour observed in the low Mach number region.

Overall, it appears that compressibility has a stabilisation effect on the transition to one of the two branches. The bifurcation point, in fact is gradually translated towards smaller viscosities for increasing values of the Mach number. In previous work by Karantonis et al.\cite{karantonis2021compressibility}, the same conclusions were drawn. An intuitive explanation can be seek in dimensional analysis. The Reynolds number, defined as the ratio between inertial and viscous forces, drives the instability. In the present geometric configuration, for sufficiently large Mach numbers, the flow is expected to undergo a more or less strong thermodynamic expansion characterised by a decrease of the gas density. Consequently, for the same value of viscosity, inertial forces will be less relevant for increasing values of the Mach number. 
Finally, it is interesting to analyse in more detail the limit of large Mach numbers. It is expected that for these values, vast regions of the domain will undergo supersonic conditions, allowing the natural development of shock waves and leading to a consistently different structure of the solution.

\begin{figure}
\centering
\begin{overpic}[trim=0 400 0 100,clip,width=0.95\columnwidth]{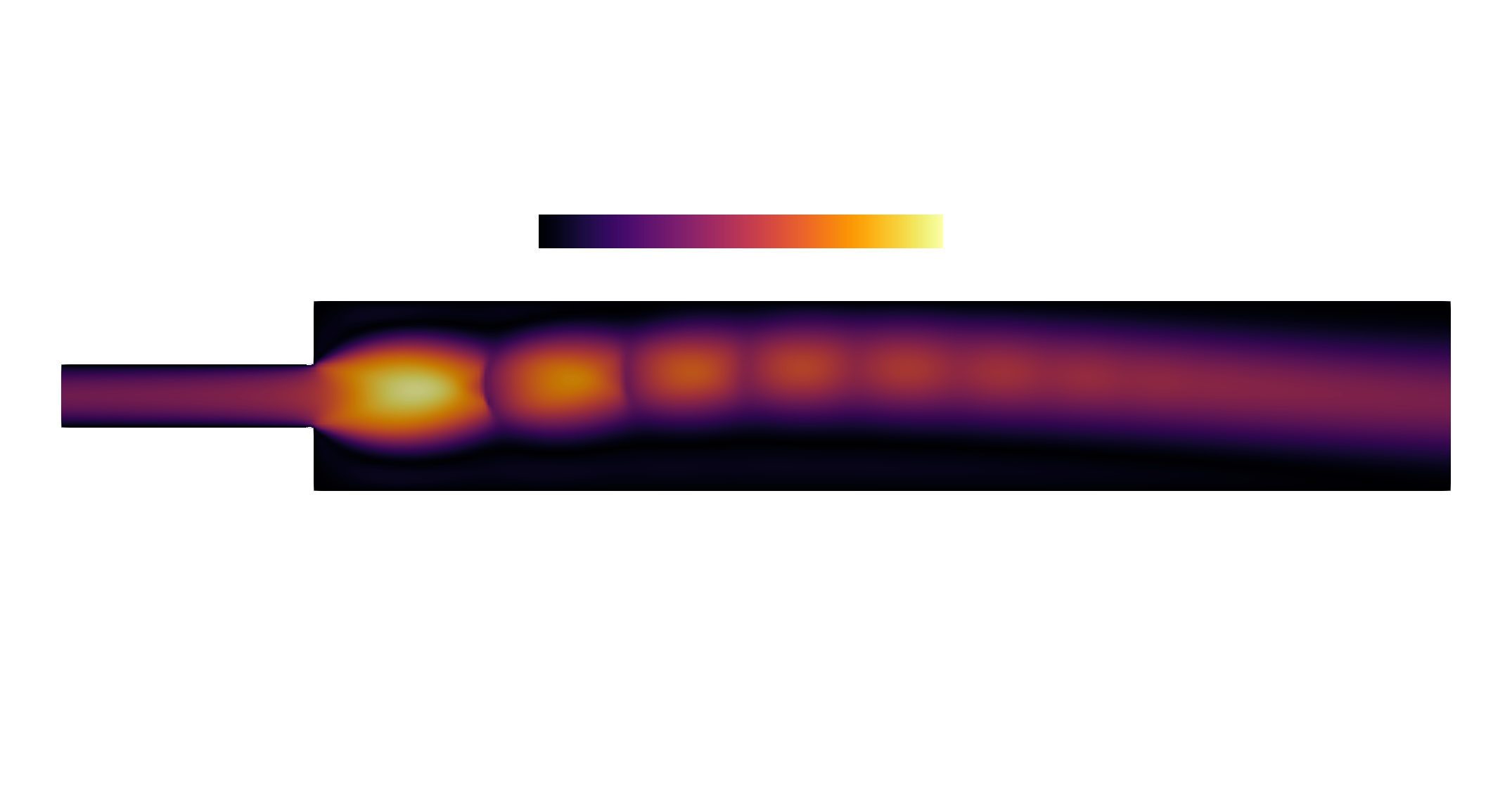}
\put (33,20.5) {$0.0$}
\put (46,22) {$\mathrm{Ma}$}
\put (59,20.5) {$2.5$}
\end{overpic}
\caption{Local Mach number for the $\mathrm{Ma}=0.9$ and $\mu=0.2$ case.}
\label{fig:ma09_mu02p}
\end{figure}

In particular, in Figure \ref{fig:ma09_mu02p} the local Mach number of the steady state solution for the largest Mach number and smallest viscosity case is shown.
First of, as expected, close to the expansion, the flow becomes supersonic. Due to this behaviour, an organised structure of concatenated sharp variations of the Mach number are observed in the same region. This structures are commonly known as \emph{Mach Diamonds} and they are common feature of exhausted nozzle flows. Of course, there is a lot of interest in the aeronautical community to accurately predict these structures which can strongly vary depending on the geometry of the problem or on the thermodynamic working conditions. The natural development of shocks and Mach diamonds allows interesting analysis in terms of reduced order modelling too.

It is well known that advected-dominated problems can be particularly challenging for classical linear reduced-order models \cite{cagniart2019model}. Consequently, the present test case allows an ideal transition between an almost incompressible set-up, where standard ROM techniques can be applied, to increasingly compressible conditions where more sophisticated approaches might be needed.
%%%%%%%%%%%%%%%%%%%%%%%%%%%%%%%%%%%%%%%%%%%%%
%%%%%%%%%%%%%%%%%%%%%%%%%%%%%%%%%%%%%%%%%%%%%
%%%%%%%%%%%%%%%%%%%%%%%%%%%%%%%%%%%%%%%%%%%%%
\subsubsection{Reduced order modelling}
%%%%%%%%%%%%%%%%%%%%%%%%%%%%%%%%%%%%%%%%%%%%%
%%%%%%%%%%%%%%%%%%%%%%%%%%%%%%%%%%%%%%%%%%%%%
%%%%%%%%%%%%%%%%%%%%%%%%%%%%%%%%%%%%%%%%%%%%% 
Two different non-intrusive techniques have been employed in this part of the work: one based on Radial Basis Function (RBF) interpolation and one based on Artificial Neural Networks (ANN).
Both methodologies are implemented within the open source library \emph{EZyRB} \cite{demo18ezyrb}. 

In particular, a smoothness parameter equal to $0.1$ has been chosen for the RBF interpolation. The ANN, instead, consists in a standard feed-forward neural network composed by $2$ layers, a first layer containing $50$ neurons and a second layer of $10$ neurons. The same architecture has been used for both pressure and velocity field. The SiLu and LogSigmoid activation functions have been used respectively for vertical velocity and pressure in order to introduce non-linearity in the network. For this preliminary test, only the upper branch has been considered.

\begin{figure}
\centering
\includegraphics[width=0.9\textwidth]{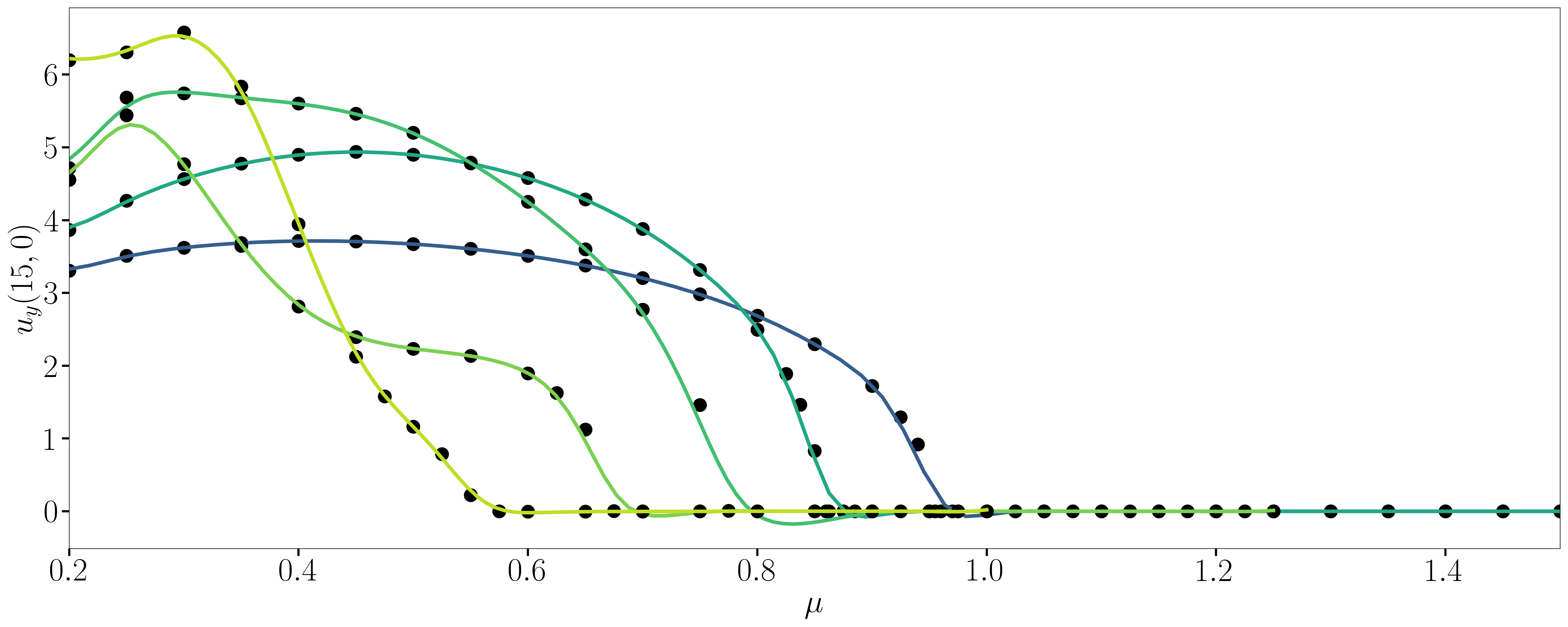}
\caption{Bifurcation diagram using RBF interpolation. Black dots indicate the full-order model snapshots. Color gradient (from blue to yellow) indicates the reduced order model for increasing values of the Mach number.}
\label{fig:diagram_RBF}
\end{figure}

\begin{figure}
\centering
\includegraphics[width=0.9\textwidth]{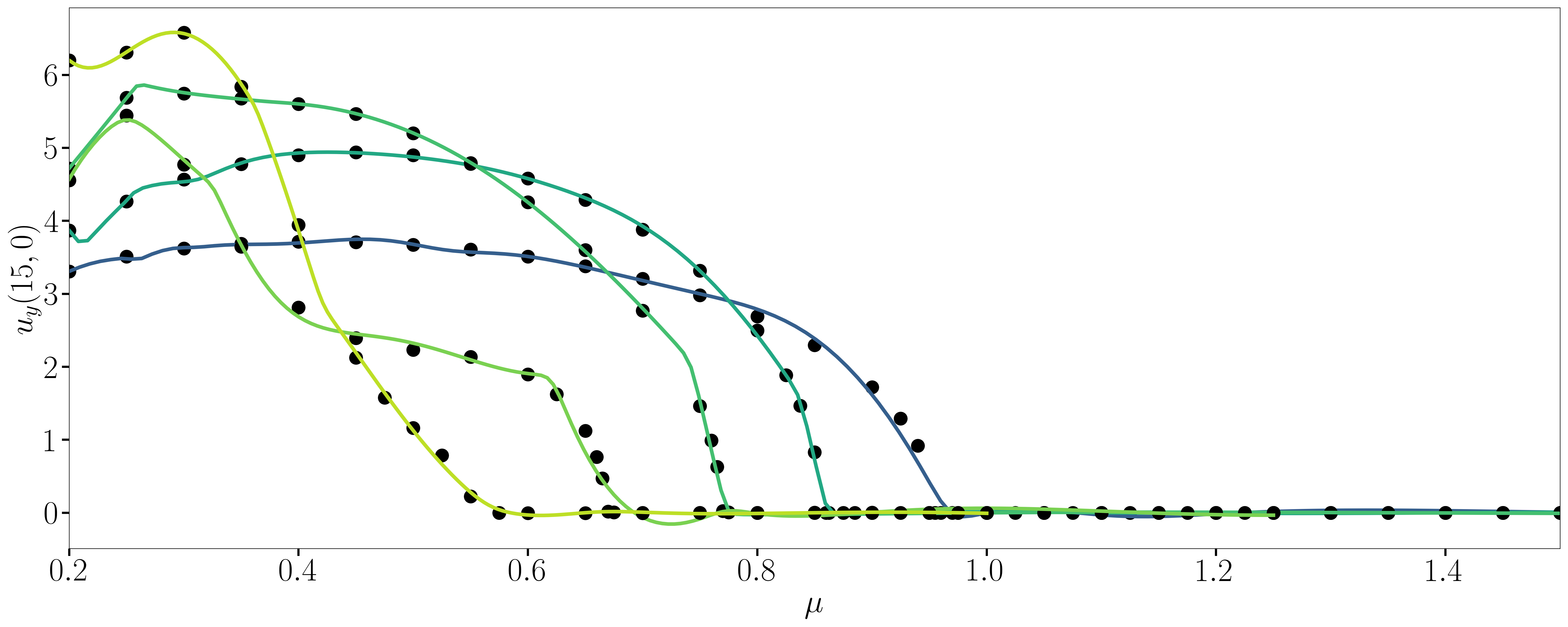}
\caption{Bifurcation diagram using ANN interpolation. Red dots indicate the full-order model snapshots. Color gradient (from blue to yellow) indicates the reduced order model for increasing values of the Mach number.}
\label{fig:diagram_ANN}
\end{figure}

\begin{figure}
\centering
\begin{overpic}[trim=0 0 400 0, clip,width=0.9\textwidth]{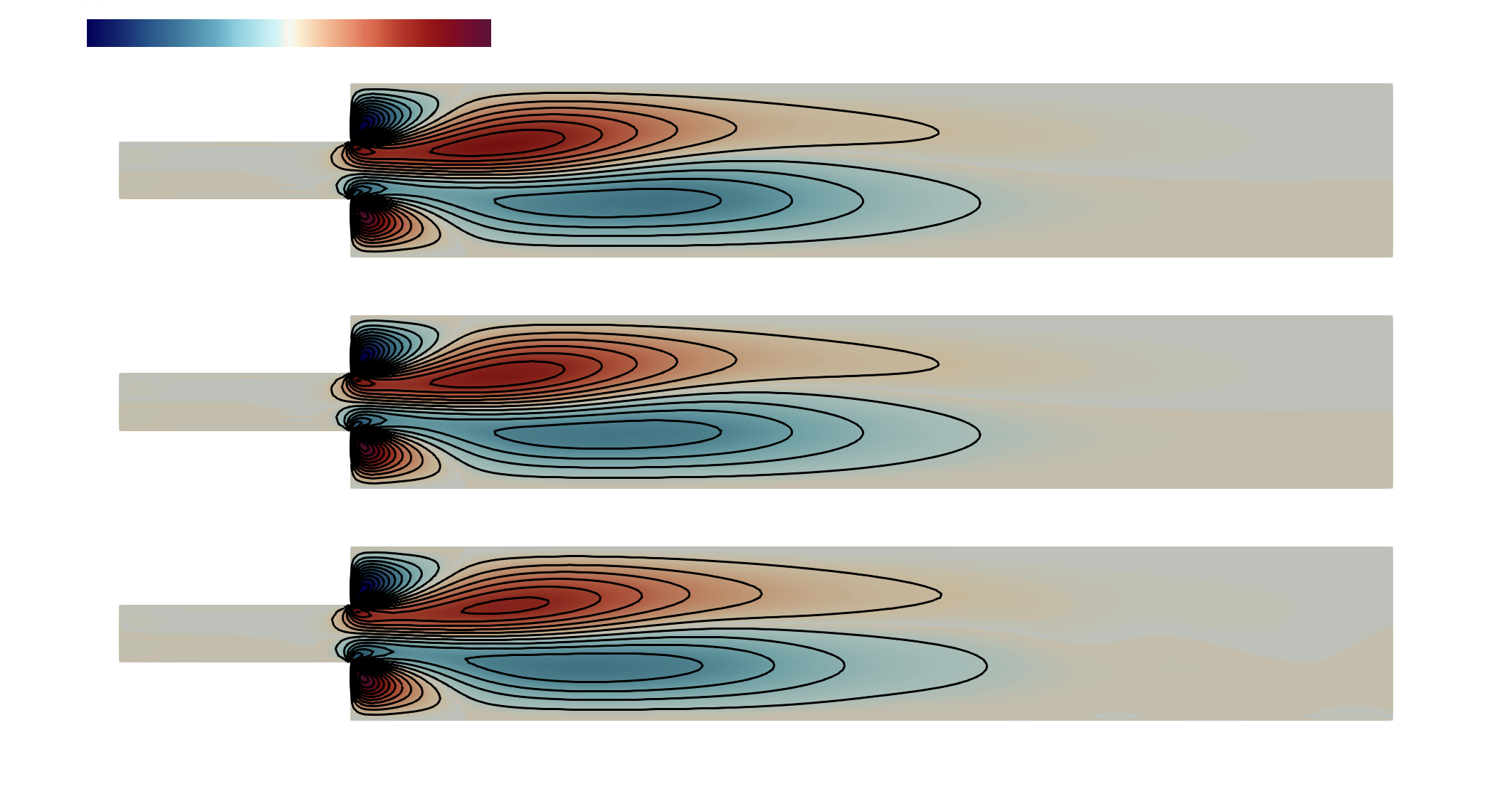}
\put (22,67.5) {$u_{y}$}
\put (2,65) {$-2.0$}
\put (36,65) {$2.0$}
\put (0,50.7) {FOM}
\put (-9,31.8) {ROM RBF}
\put (-9,13.0) {ROM ANN}
\end{overpic}
\caption{Comparison of the wall-normal velocity field between full-order model and different reduced-order models (for $\mu=0.95$). Top, full order model; middle, reduced-order model using RBF interpolation; bottom, reduced-order model using ANN interpolation.}
\label{fig:comparison_mu095v}
\end{figure}

The bifurcation diagrams reconstructed using the ROM methodologies are reported in Figures \ref{fig:diagram_RBF} (RBF) and \ref{fig:diagram_ANN} (ANN).
In a more qualitative way, using a leave-one-out strategy on the parameter space, full-order and reduced-order models have been compared. In particular, the vertical velocity field is compared for the case of $\mu=0.95$ in Figure \ref{fig:comparison_mu095v}. Observing the bifurcation diagram and the visual comparison, we can notice a good agreement between full-order model and reduced order methodologies. Keep in mind that the total number of snapshots is rather small, in particular in proximity of the bifurcation point.
\begin{figure}
\centering
\begin{overpic}[trim=0 0 400 0, clip,width=0.9\textwidth]{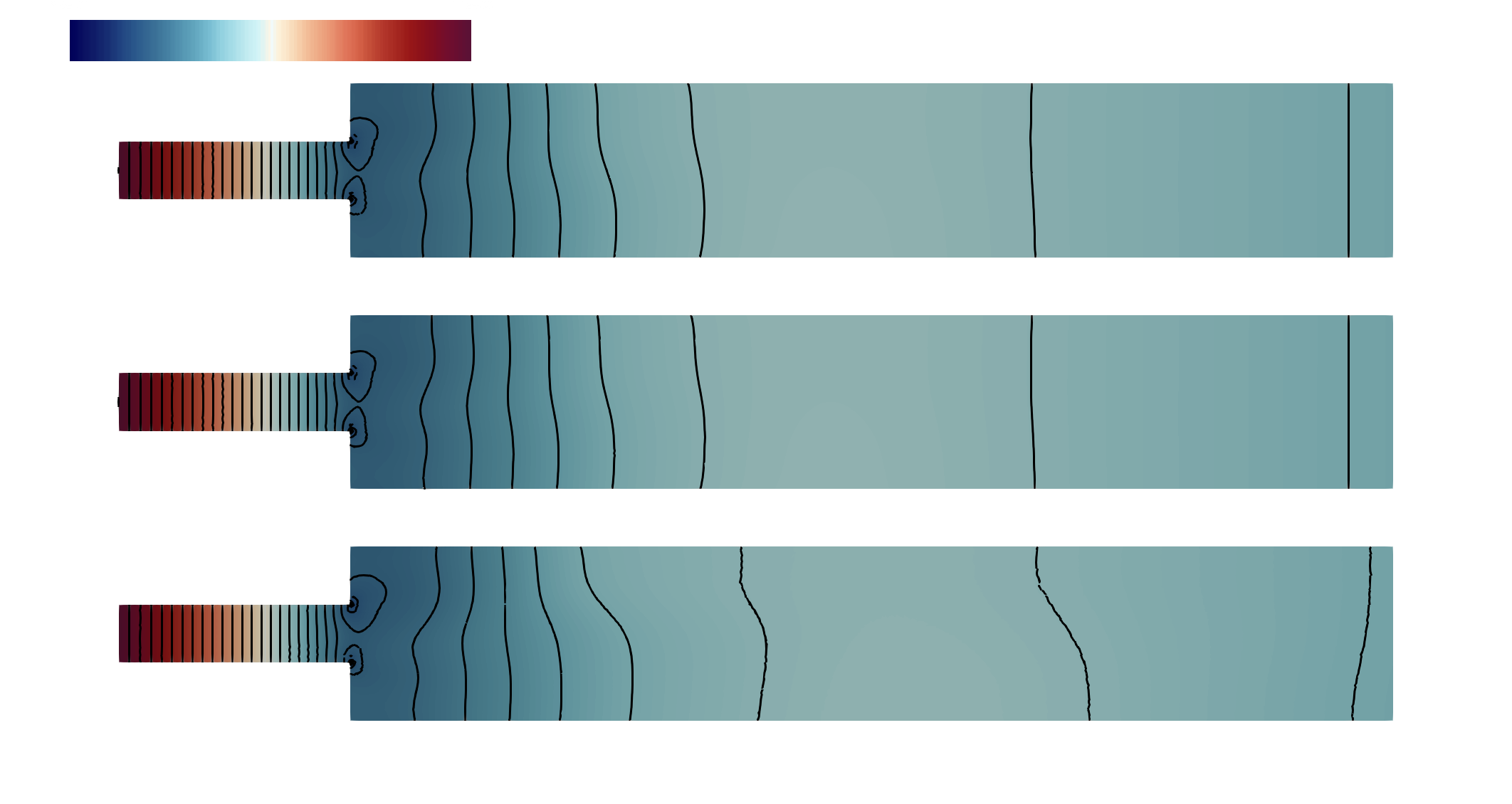}
\put (18.5,67.5) {$p/p_{\textrm{out}}$}
\put (2,65) {$0.95$}
\put (35,65) {$1.10$}
\put (0,50.7) {FOM}
\put (-9,31.8) {ROM RBF}
\put (-9,13.0) {ROM ANN}
\end{overpic}
\caption{Comparison of the pressure field between full-order model and different reduced-order models (for $\mu=0.95$ and $\textrm{Ma}=0.3$). Top, full order model; middle, reduced-order model using RBF interpolation; bottom, reduced-order model using ANN interpolation.}
\label{fig:comparison_mu095p}
\end{figure}
A similar comparison for the pressure field in shown in Figure \ref{fig:comparison_mu095p}.

\begin{figure}
\centering
\subfigure[Vertical velocity.]{\includegraphics[width=0.48\textwidth]{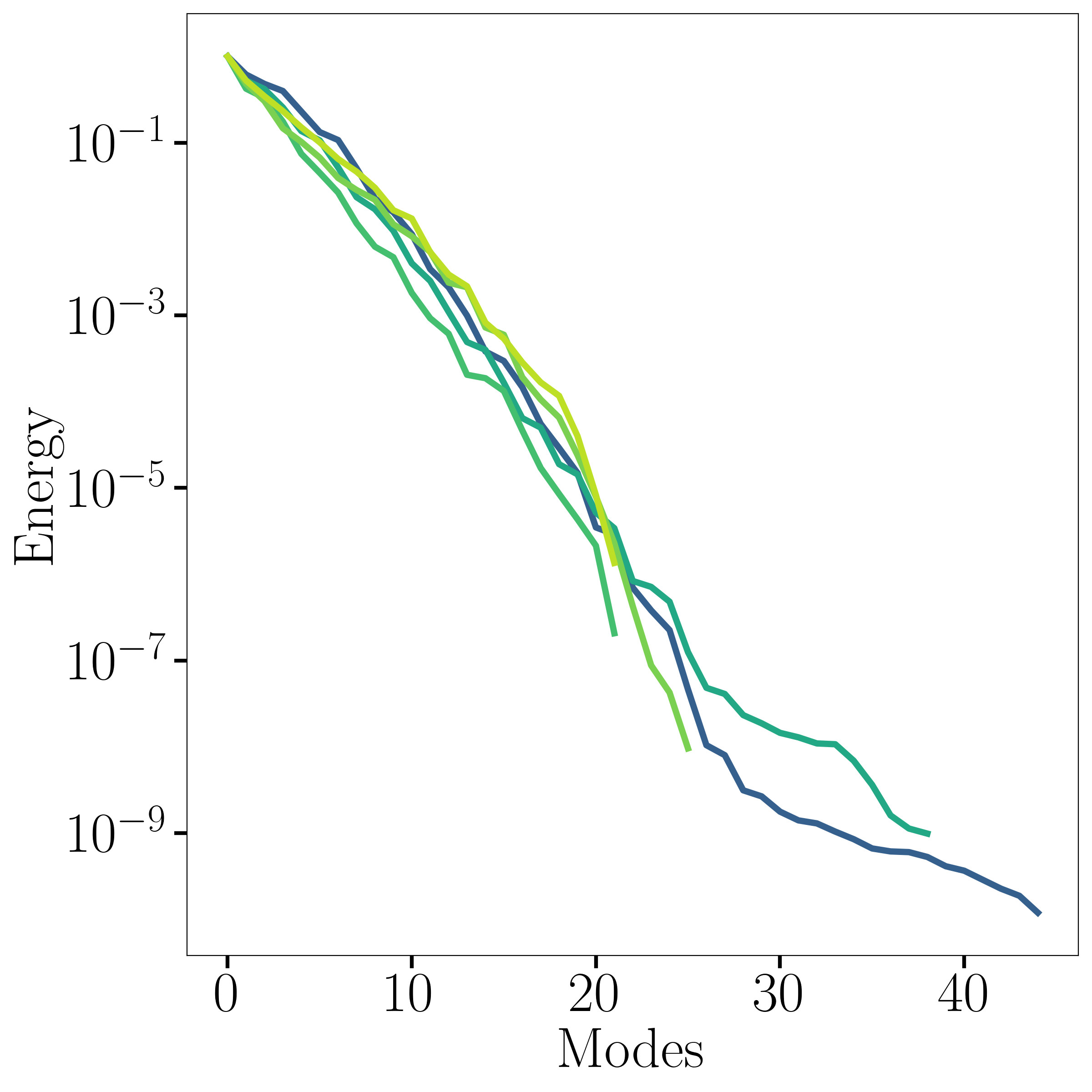}\label{fig:singular_values_v}}
\subfigure[Pressure.]{\includegraphics[width=0.48\textwidth]{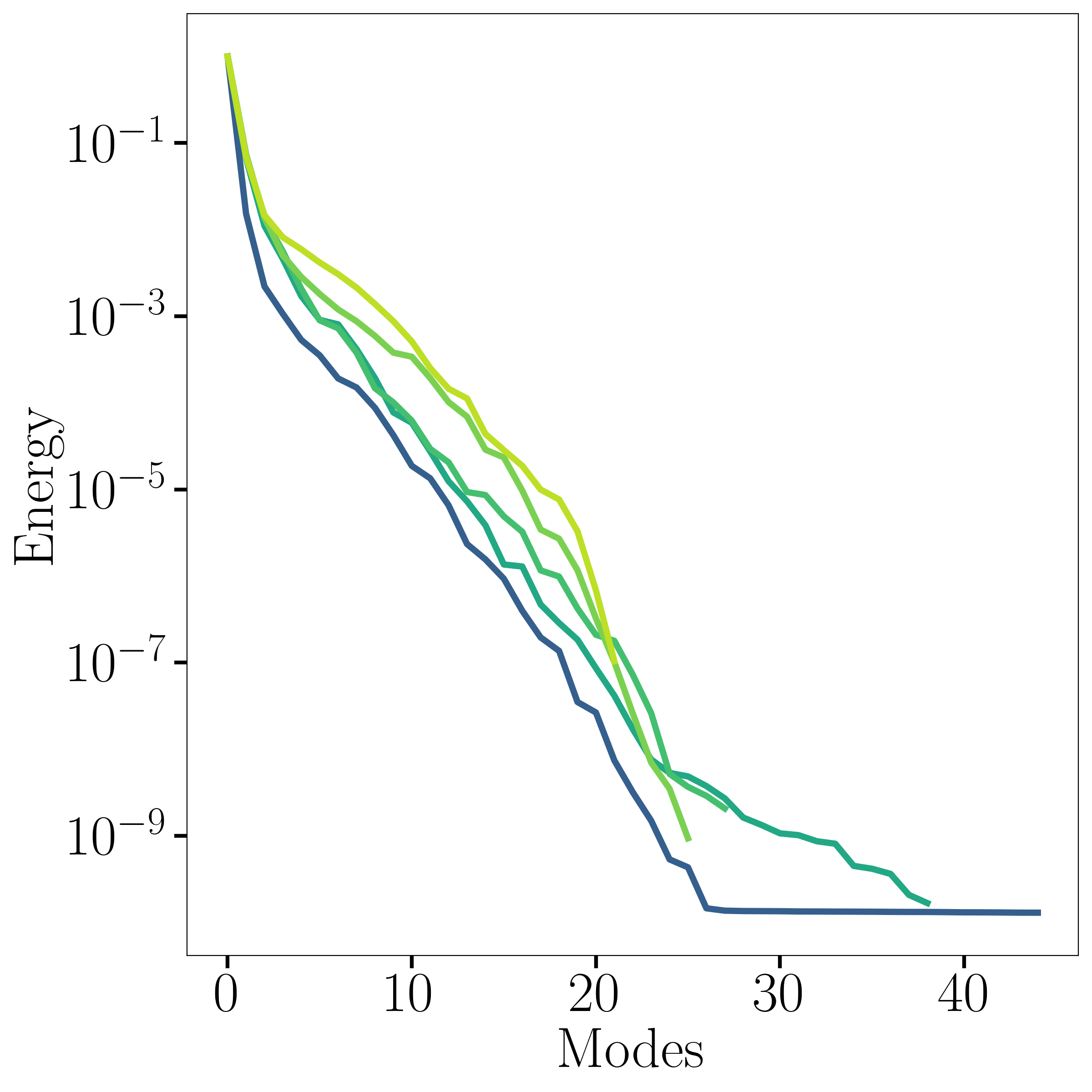}\label{fig:singular_values_p}}
\caption{Normalised energy of the singular values from the POD of the matrix of snapshots for vertical velocity (left) and pressure (right) fields. Color gradient (from blue to yellow) indicates increasing values of the Mach number.}
\label{fig:singular_values}
\end{figure}

\begin{figure}
\centering
\includegraphics[ width=0.9\textwidth]{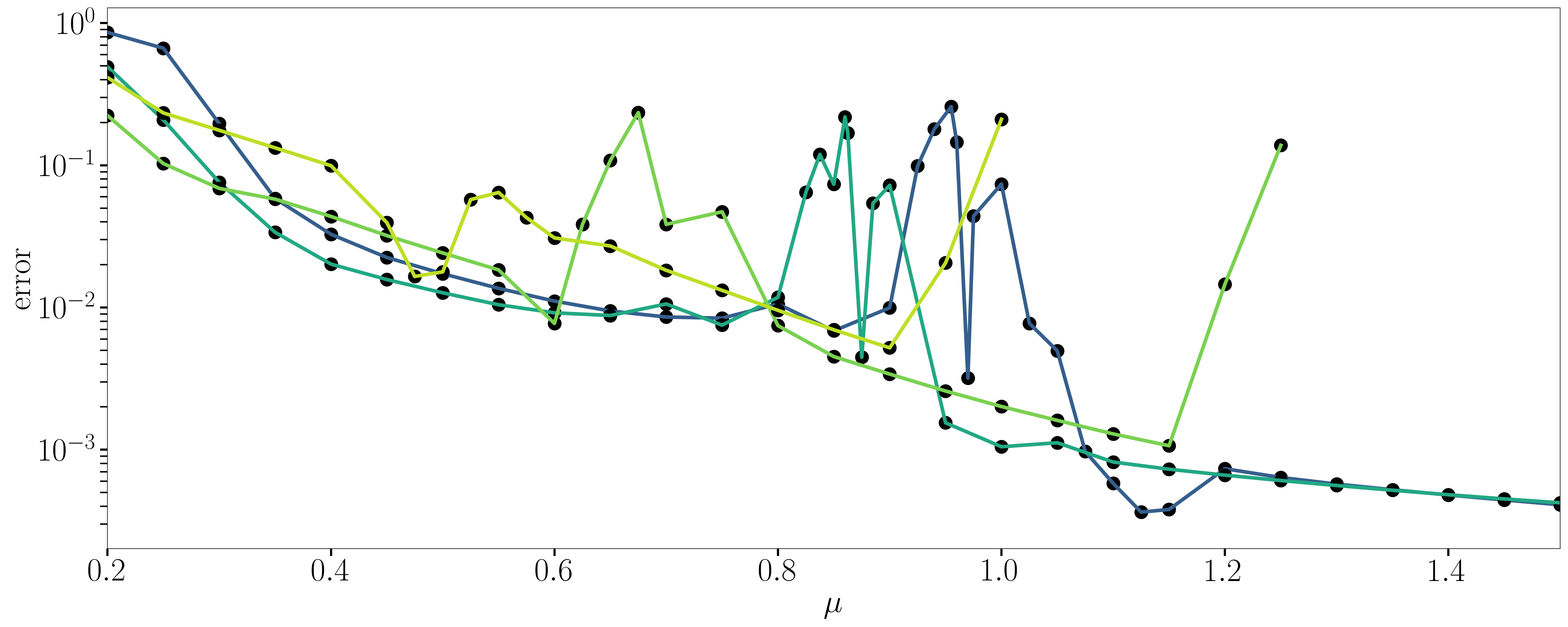}
\caption{Leave-one-out relative error for the vertical component of the velocity field using the RBF interpolation. Color gradient (from blue to yellow) indicates increasing values of the Mach number.}
\label{fig:loo_error_RBFv}
\end{figure}

\begin{figure}
\centering
\includegraphics[ width=0.9\textwidth]{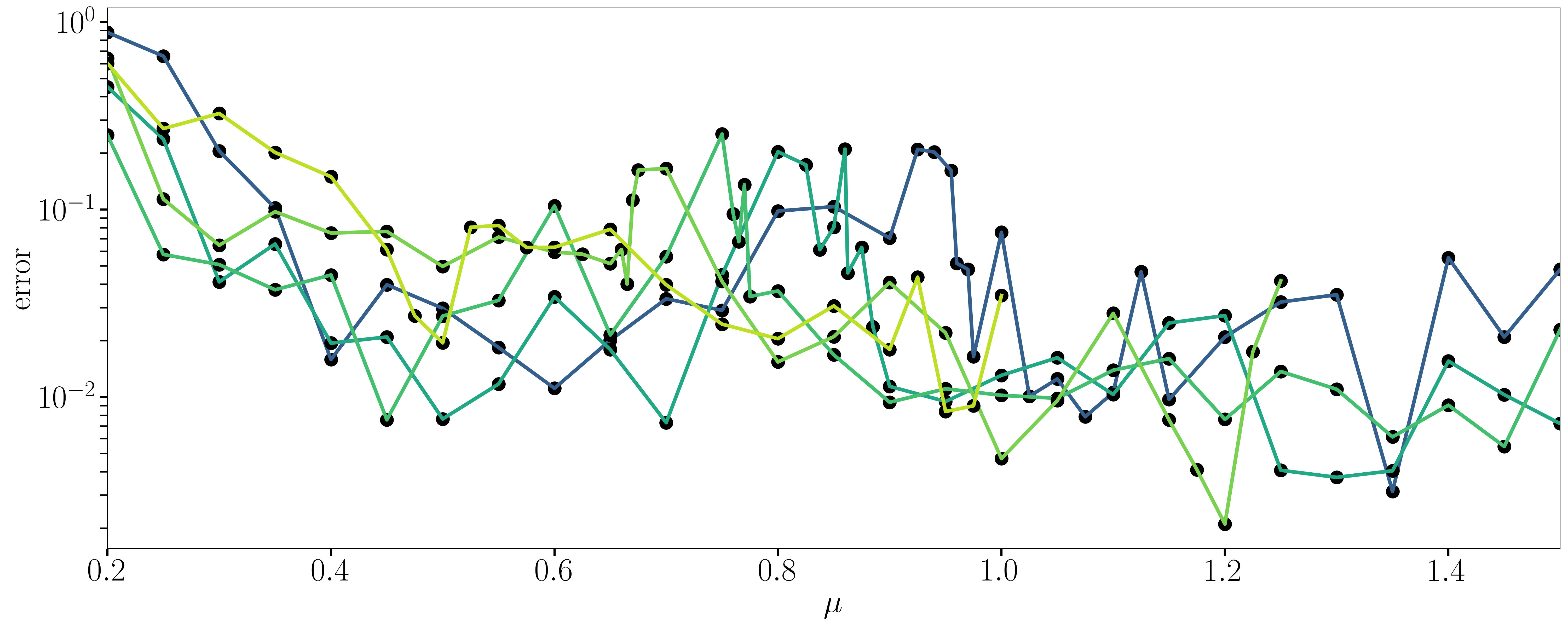}
\caption{Leave-one-out relative error for the vertical component of the velocity field using the ANN interpolation. Color gradient (from blue to yellow) indicates increasing values of the Mach number.}
\label{fig:loo_error_ANNv}
\end{figure}

\begin{figure}
\centering
\includegraphics[ width=0.9\textwidth]{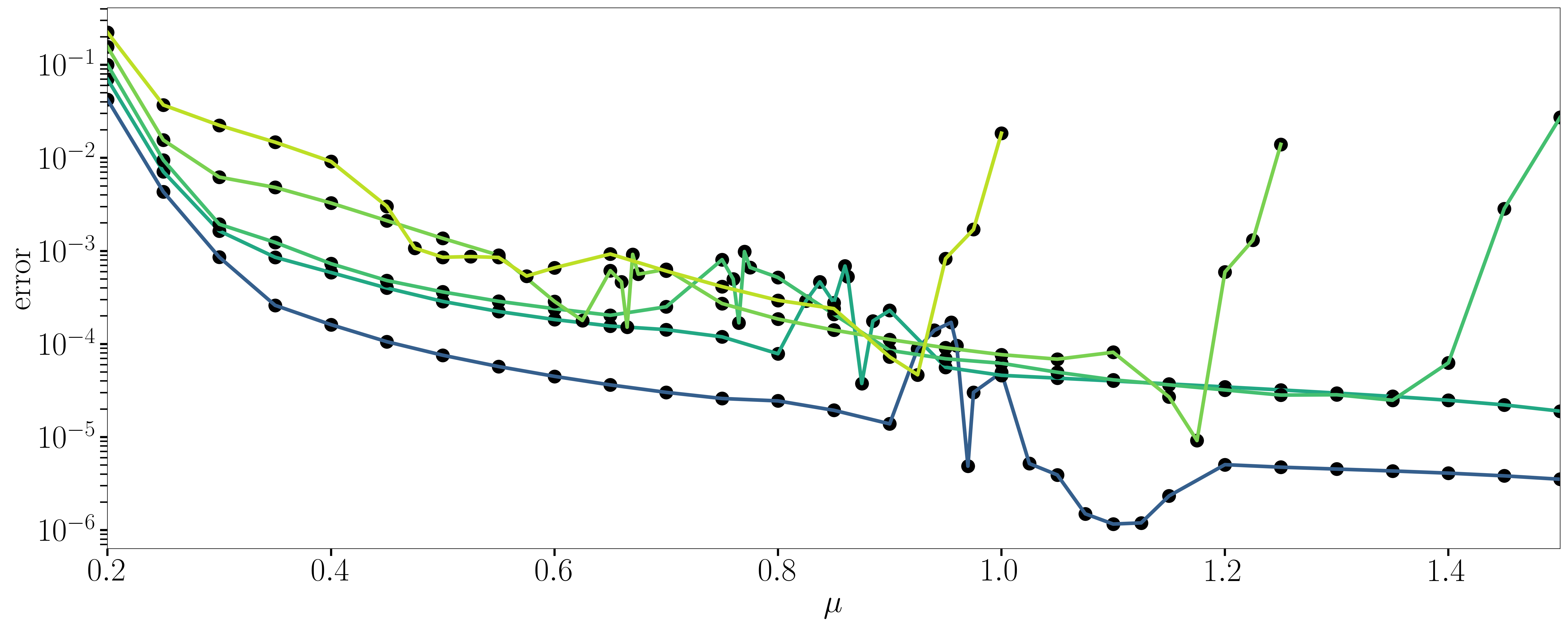}
\caption{Leave-one-out relative error for the pressure using the RBF interpolation. Color gradient (from blue to yellow) indicates increasing values of the Mach number.}
\label{fig:loo_error_RBFp}
\end{figure}

\begin{figure}
\centering
\includegraphics[ width=0.9\textwidth]{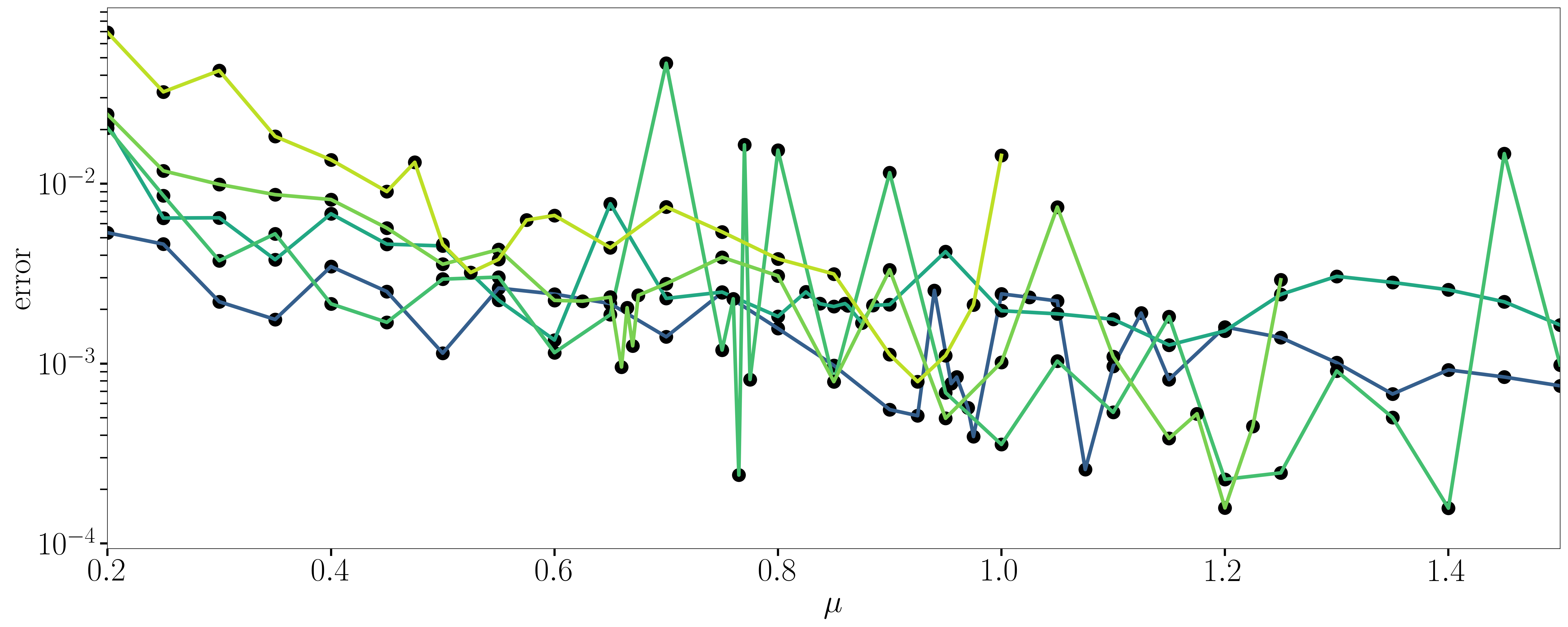}
\caption{Leave-one-out relative error for the pressure using the ANN interpolation. Color gradient (from blue to yellow) indicates increasing values of the Mach number.}
\label{fig:loo_error_ANNp}
\end{figure}

\begin{figure}
\centering
\includegraphics[width=0.9\textwidth]{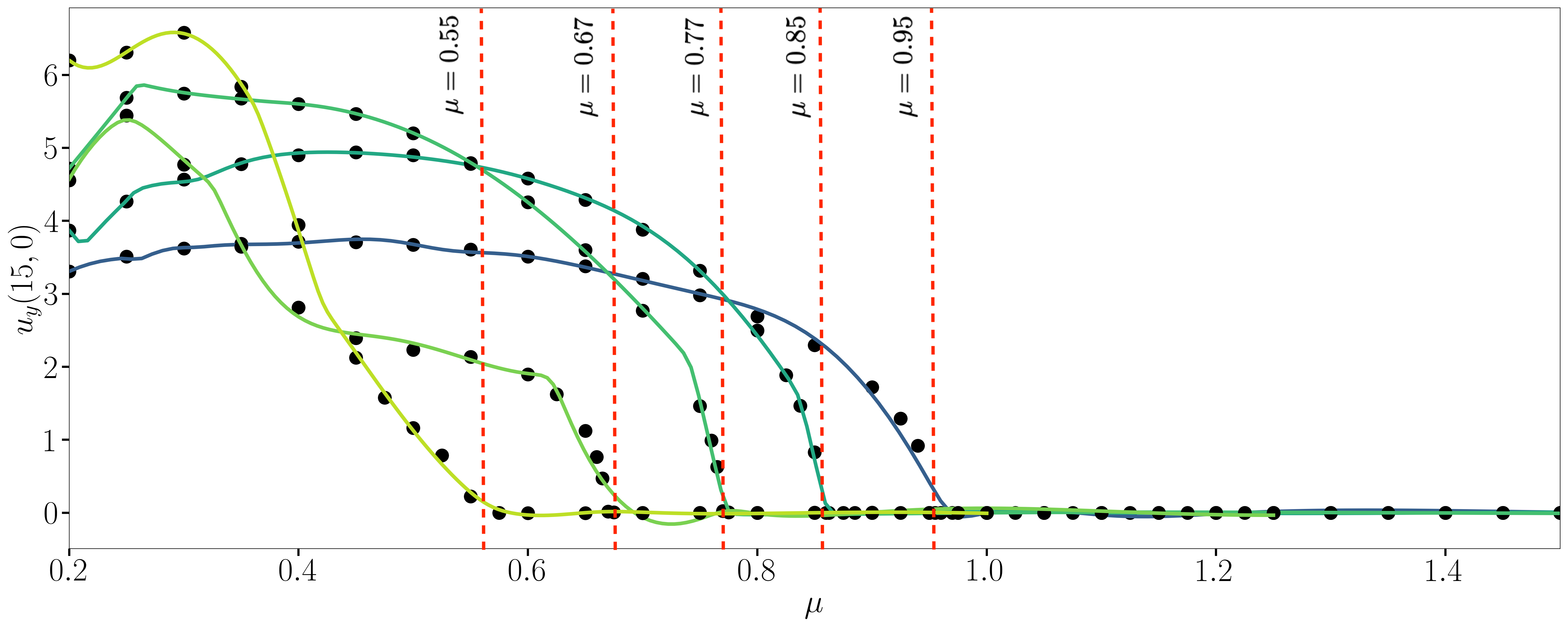}
\caption{Bifurcation diagram using ANN interpolation. Red dots indicate the full-order model snapshots. The solid blue line represents the reconstructed diagram using the ROM approach.}
\label{fig:diagram_ANN2}
\end{figure}

\begin{figure}
\centering
\begin{overpic}[trim=0 0 400 0, clip,width=0.9\textwidth]{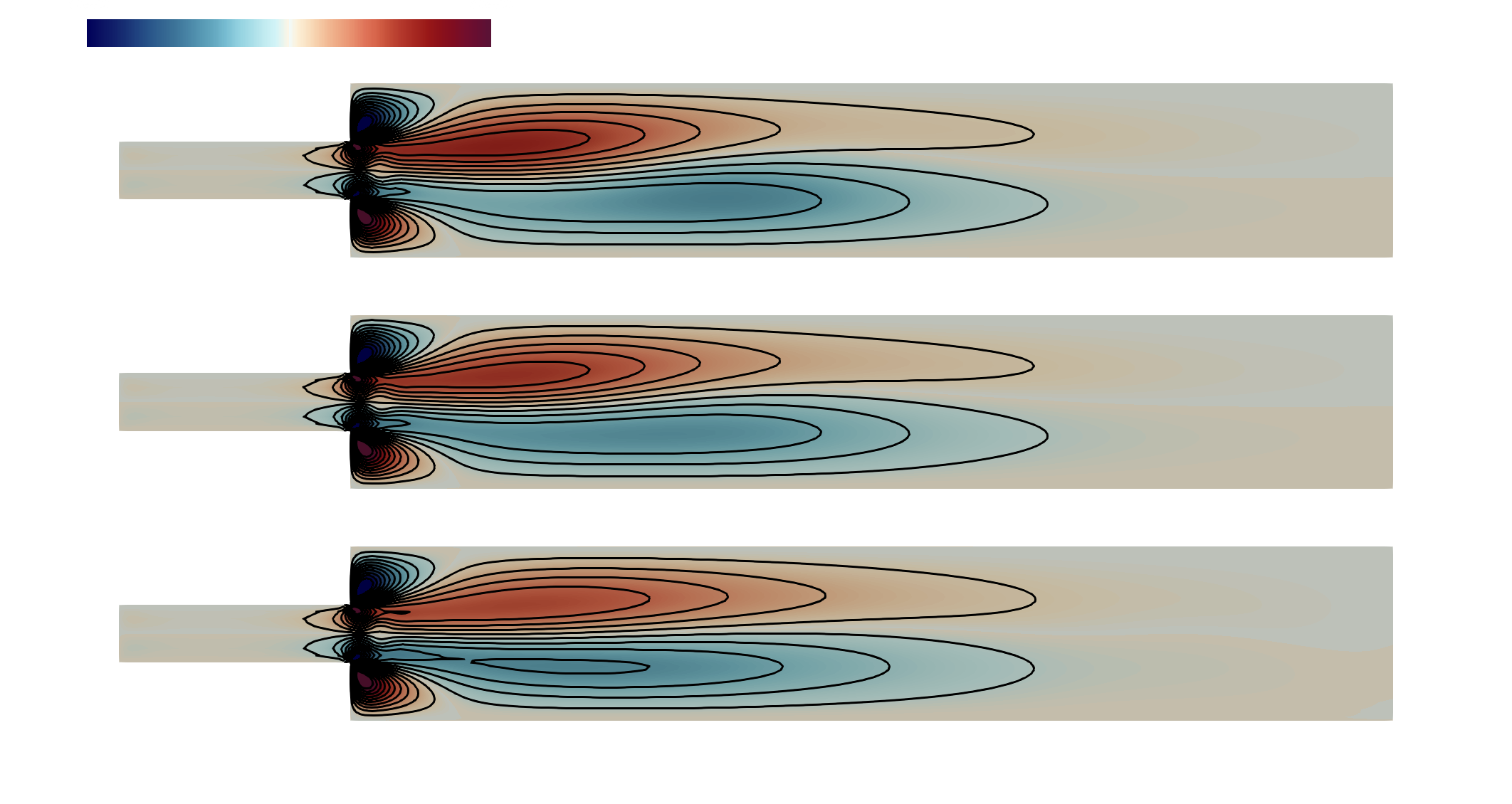}
\put (22,67.5) {$u_{y}$}
\put (2,65) {$-3.0$}
\put (40,65) {$3.0$}
\put (0,50.7) {FOM}
\put (-9,31.8) {ROM RBF}
\put (-9,13.0) {ROM ANN}
\end{overpic}
\caption{Comparison of the wall-normal velocity field between full-order model and different reduced-order models (for $\mu=0.855$  and $\textrm{Ma}=0.6$). Top, full order model; middle, reduced-order model using RBF interpolation; bottom, reduced-order model using ANN interpolation.}
\label{fig:comparison_mu0855v}
\end{figure}

\begin{figure}
\centering
\begin{overpic}[trim=0 0 400 0, clip,width=0.9\textwidth]{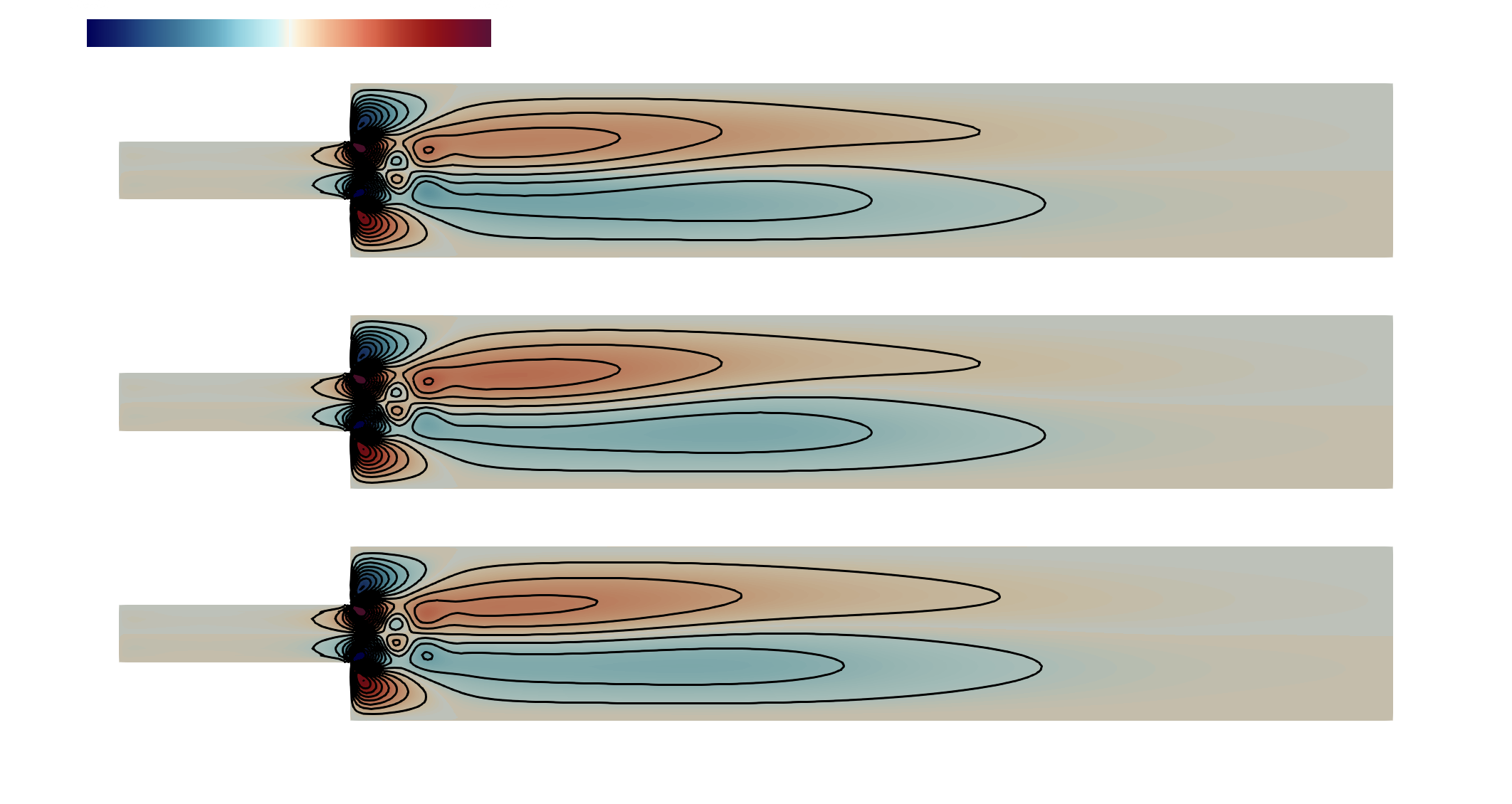}
\put (22,67.5) {$u_{y}$}
\put (2,65) {$-5.0$}
\put (37,65) {$5.0$}
\put (0,50.7) {FOM}
\put (-9,31.8) {ROM RBF}
\put (-9,13.0) {ROM ANN}
\end{overpic}
\caption{Comparison of the wall-normal velocity field between full-order model and different reduced-order models (for $\mu=0.77$  and $\textrm{Ma}=0.7$). Top, full order model; middle, reduced-order model using RBF interpolation; bottom, reduced-order model using ANN interpolation.}
\label{fig:comparison_mu077v}
\end{figure}

\begin{figure}
\centering
\begin{overpic}[trim=0 0 400 0, clip,width=0.9\textwidth]{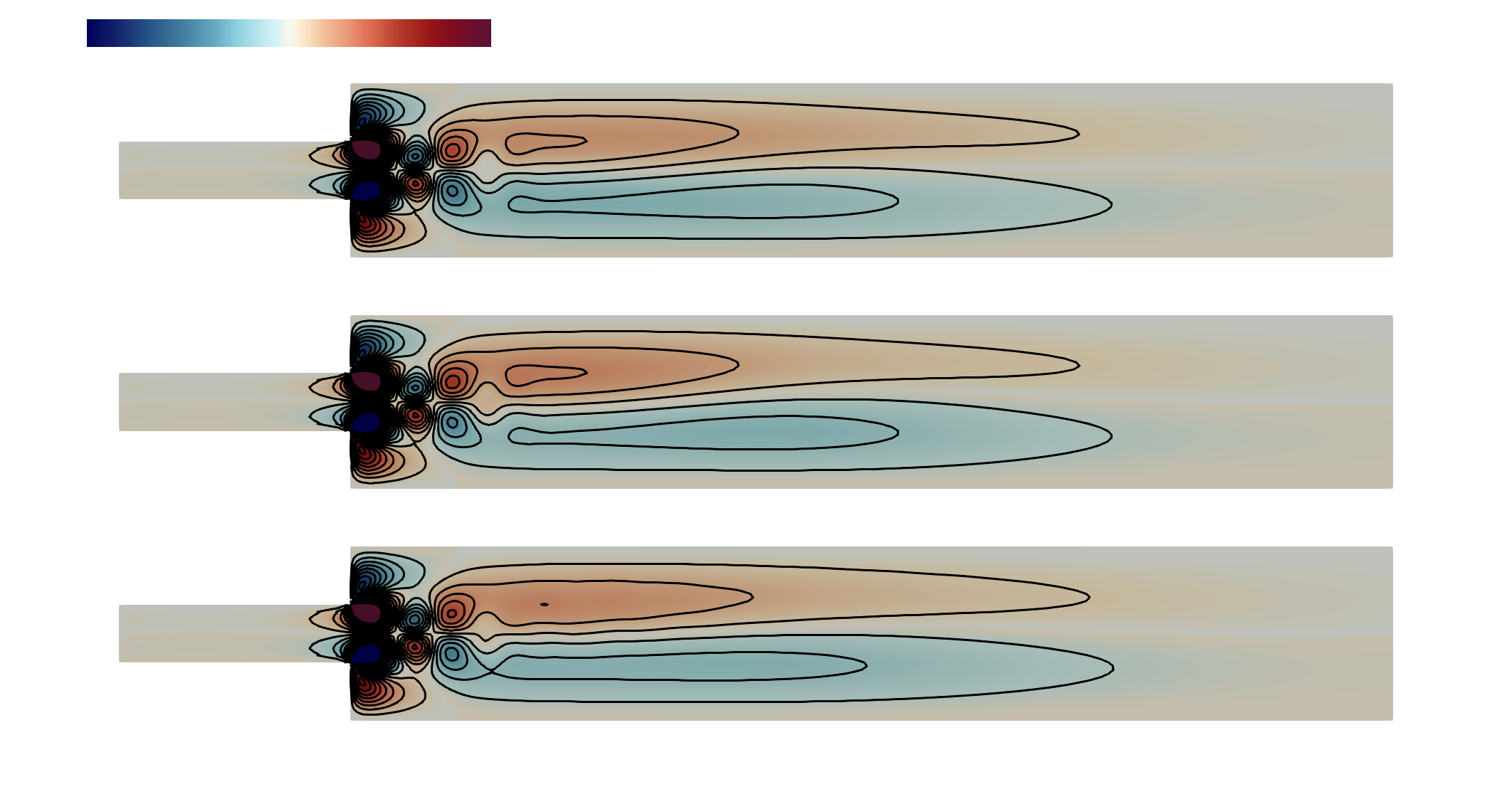}
\put (22,67.5) {$u_{y}$}
\put (2,65) {$-5.0$}
\put (36,65) {$5.0$}
\put (0,50.7) {FOM}
\put (-9,31.8) {ROM RBF}
\put (-9,13.0) {ROM ANN}
\end{overpic}
\caption{Comparison of the wall-normal velocity field between full-order model and different reduced-order models (for $\mu=0.67$  and $\textrm{Ma}=0.8$). Top, full order model; middle, reduced-order model using RBF interpolation; bottom, reduced-order model using ANN interpolation.}
\label{fig:comparison_mu067v}
\end{figure}

\begin{figure}
\centering
\begin{overpic}[trim=0 0 400 0, clip,width=0.9\textwidth]{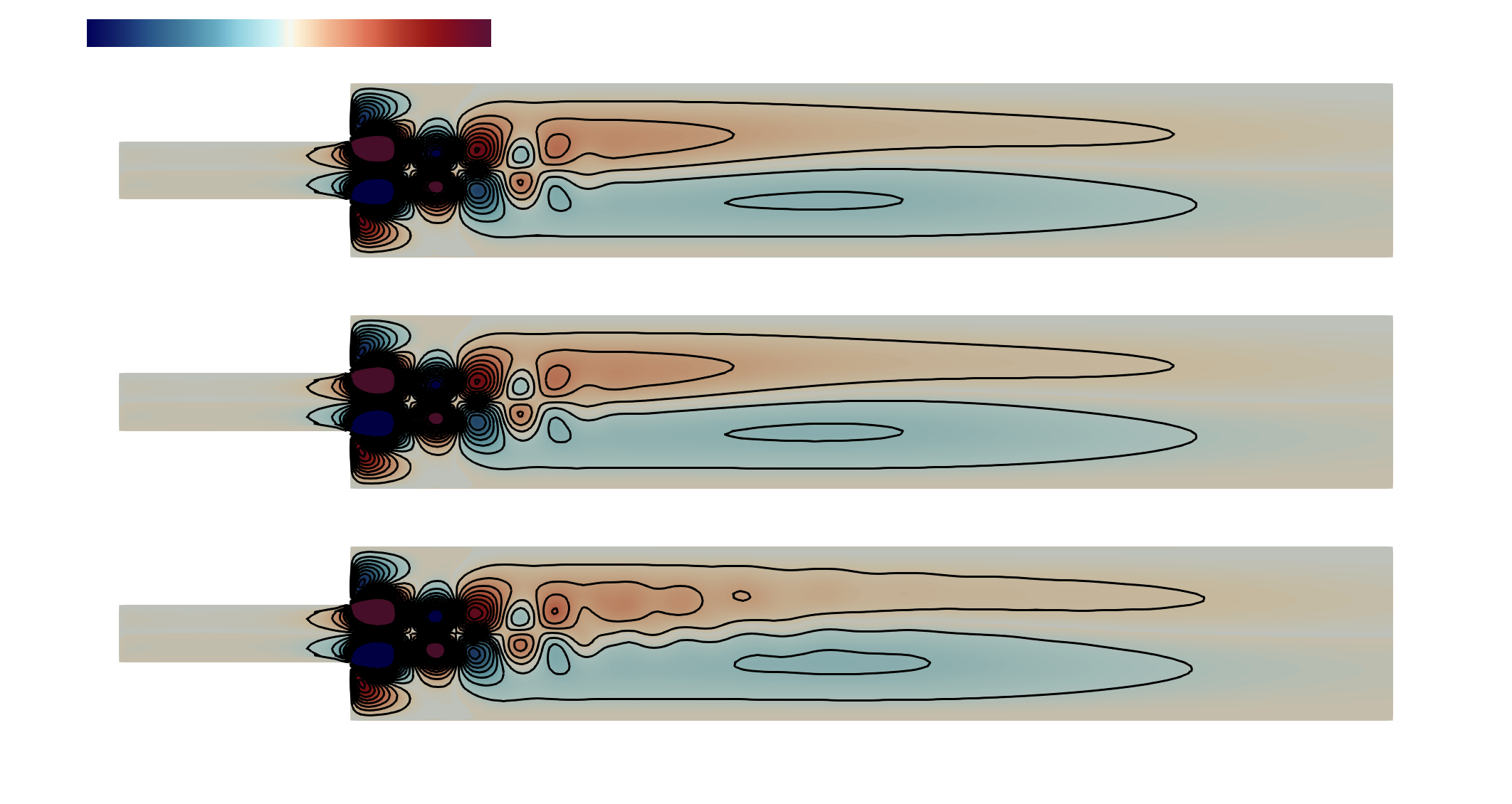}
\put (22,67.5) {$u_{y}$}
\put (2,65) {$-5.0$}
\put (37,65) {$5.0$}
\put (0,50.7) {FOM}
\put (-9,31.8) {ROM RBF}
\put (-9,13.0) {ROM ANN}
\end{overpic}
\caption{Comparison of the wall-normal velocity field between full-order model and different reduced-order models (for $\mu=0.56$  and $\textrm{Ma}=0.9$). Top, full order model; middle, reduced-order model using RBF interpolation; bottom, reduced-order model using ANN interpolation.}
\label{fig:comparison_mu056v}
\end{figure}

In order to better understand the influence of compressibility on the reduced order model, the singular values of the most relevant quantities have been computed. The decay of the singular values for vertical velocity and pressure at different Mach numbers is shown in Figure \ref{fig:singular_values}. 
A few conclusions can be drawn from this representation. First, in the lowest Mach number case, the decay of the modes for the vertical component of the velocity field is quite steep up to approximately $25$ modes. After that, it slows down significantly. For the pressure field, at the same Mach number, instead, after $25$ modes the energy reaches a plateau. 
In terms of the influence of compressibility, it appears that increasing the Mach number does not change the decay of the modes for the vertical velocity. Instead, a slower decay is observed in the pressure field for increasing values of the Mach number. This is probably due to the more central role played by internal energy (i.e. thermodynamic pressure) on the overall dynamics of the system. More complex behaviours of the pressure field are consequently expected, leading to a stiffer interpolation problem on the ROM side.
It is interesting to highlight that the dimensional reduction step is the same for both RBF and ANN approaches as they both rely on the POD of the same snapshots. The two methodologies differ in the approximation/interpolation step. Consequently, the previous discussion applies to both approaches. 
A second comparison between the two methodologies consists in comparing the leave-one-out errors for vertical velocity and pressure fields at different Mach numbers. These quantities are shown for the vertical velocity field in Figures \ref{fig:loo_error_RBFv} and \ref{fig:loo_error_ANNv} for RBF and ANN respectively. The same comparisons for the pressure field are represented in Figures \ref{fig:loo_error_RBFp} and \ref{fig:loo_error_ANNp}.
Rather clear trends can be identified in the leave-one-out error of the pressure field (Figure \ref{fig:loo_error_RBFp}). First, the error tends to increase while approaching the extrema of the parameter space, as normally expected. More interestingly, similarly to the corresponding plot for vertical velocity, the errors are larger for smaller values of viscosity. In this region, in fact, the hyperbolic nature of the Navier-Stokes system prevails, making the approximation between snapshots more challenging. In particular, for the pressure field, for increasing values of the Mach number, this tendency is more and more pronounced. The explanation is somehow similar: for larger values of the Mach number and small values of viscosity, the system is more prone to the development of shock waves, which are known to be a particularly challenging feature to be approximated by reduced order modelling \cite{nair2019transported}. The behaviour of the leave-one-out errors for the ANN is characterised by less evident patterns. In Figure \ref{fig:loo_error_ANNp}, it is still noticeable the increase of the errors approaching smaller and smaller values of viscosity. Along the same lines of Figure \ref{fig:loo_error_RBFp}, it can also be noticed the increase of the approximation errors for larger values of the Mach number in the low viscosity region.
Finally, the typical peaks of the leave-one-out errors in proximity of the bifurcation are more evident in Figures \ref{fig:loo_error_RBFv} and \ref{fig:loo_error_ANNv} rather than in the pressure field.
A comparison between full and reduced order models in proximity of the bifurcation have been considered and shown in the following figures. In particular, the critical values of viscosity are highlighted as dashed red lines in Figure \ref{fig:diagram_ANN2}. The same values are also listed in table \ref{tab:muROM}.
\begin{table}
\centering
\begin{tabular}{|c|c|c|c|c|c|}
\hline
$\Ma$ & $0.3$ & $0.6$ & $0.7$ & $0.8$ & $0.9$  \\
\hline
$\mu_{\textrm{ROM}}$ & $0.95$ & $0.855$ & $0.77$ & $0.67$  & $0.55$ \\
\hline
\end{tabular}
\caption{}
\label{tab:muROM}
\end{table}
It can be noticed that both ROM strategies are able to predict the general behaviour of the bifurcation; the chosen values of viscosity are linked to mildly asymmetric solutions. It is evident in the representation of the vertical velocity field in Figures \ref{fig:comparison_mu0855v}, \ref{fig:comparison_mu077v}, \ref{fig:comparison_mu067v} and \ref{fig:comparison_mu056v}. It can also be noticed, in all the figures, the increasing complexity of the flow field, both in terms of velocity and pressure field, for larger values of the Mach number. For example, this is observable by comparing the vertical component of the velocity field at $\textrm{Ma}=0.3$ (Figure \ref{fig:comparison_mu095v}) and $\textrm{Ma}=0.9$ (Figure \ref{fig:comparison_mu056v}).
\begin{figure}
\centering
\begin{overpic}[trim=0 0 400 0, clip,width=0.9\textwidth]{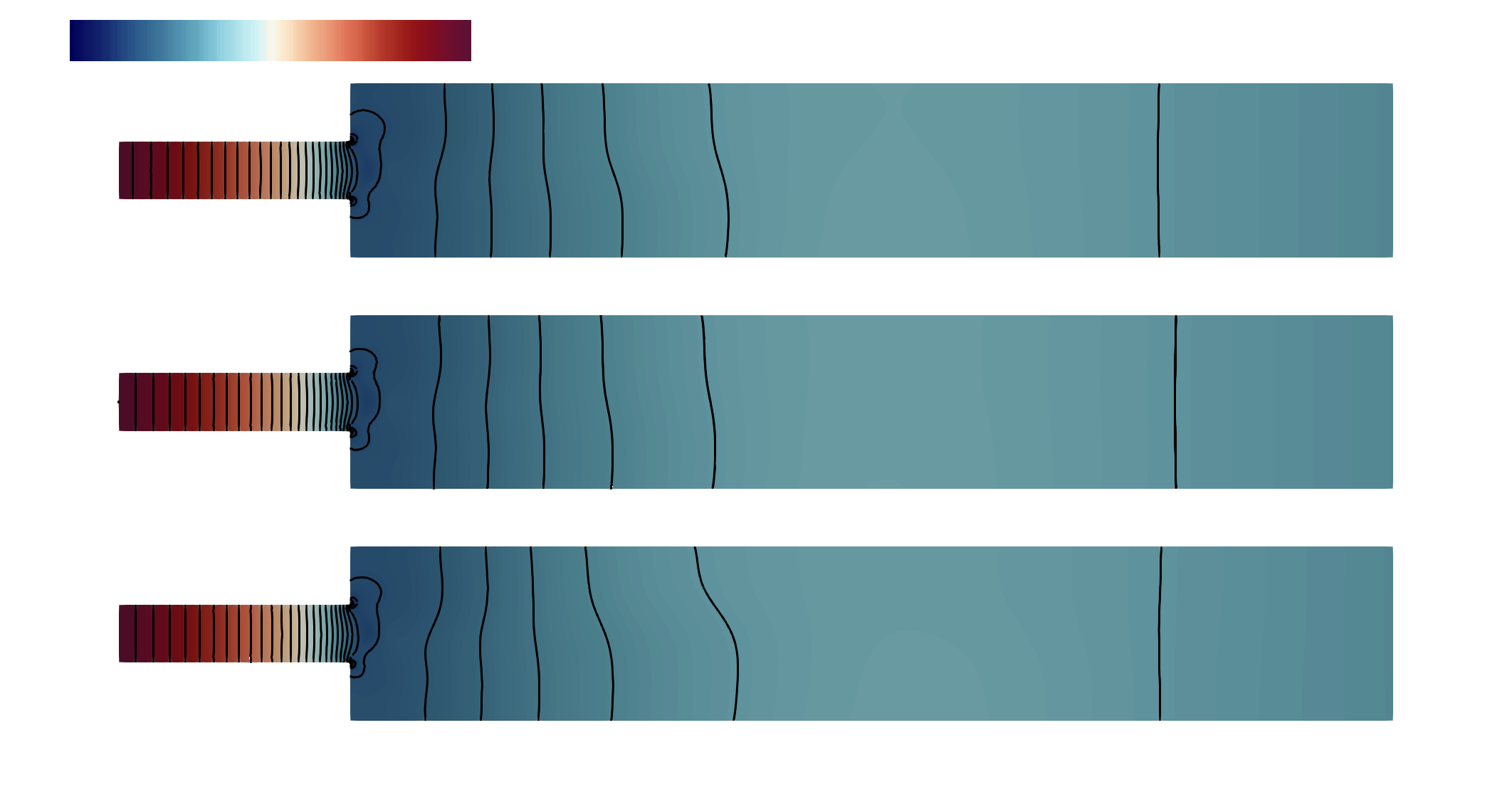}
\put (18.5,67.5) {$p/p_{\textbf{out}}$}
\put (2,65) {$0.76$}
\put (35,65) {$2.13$}
\put (0,50.7) {FOM}
\put (-9,31.8) {ROM RBF}
\put (-9,13.0) {ROM ANN}
\end{overpic}
\caption{Comparison of the pressure field between full-order model and different reduced-order models (for $\mu=0.855$ and $\textrm{Ma}=0.6$). Top, full order model; middle, reduced-order model using RBF interpolation; bottom, reduced-order model using ANN interpolation.}
\label{fig:comparison_mu0855p}
\end{figure}
\begin{figure}
\centering
\begin{overpic}[trim=0 0 400 0, clip,width=0.9\textwidth]{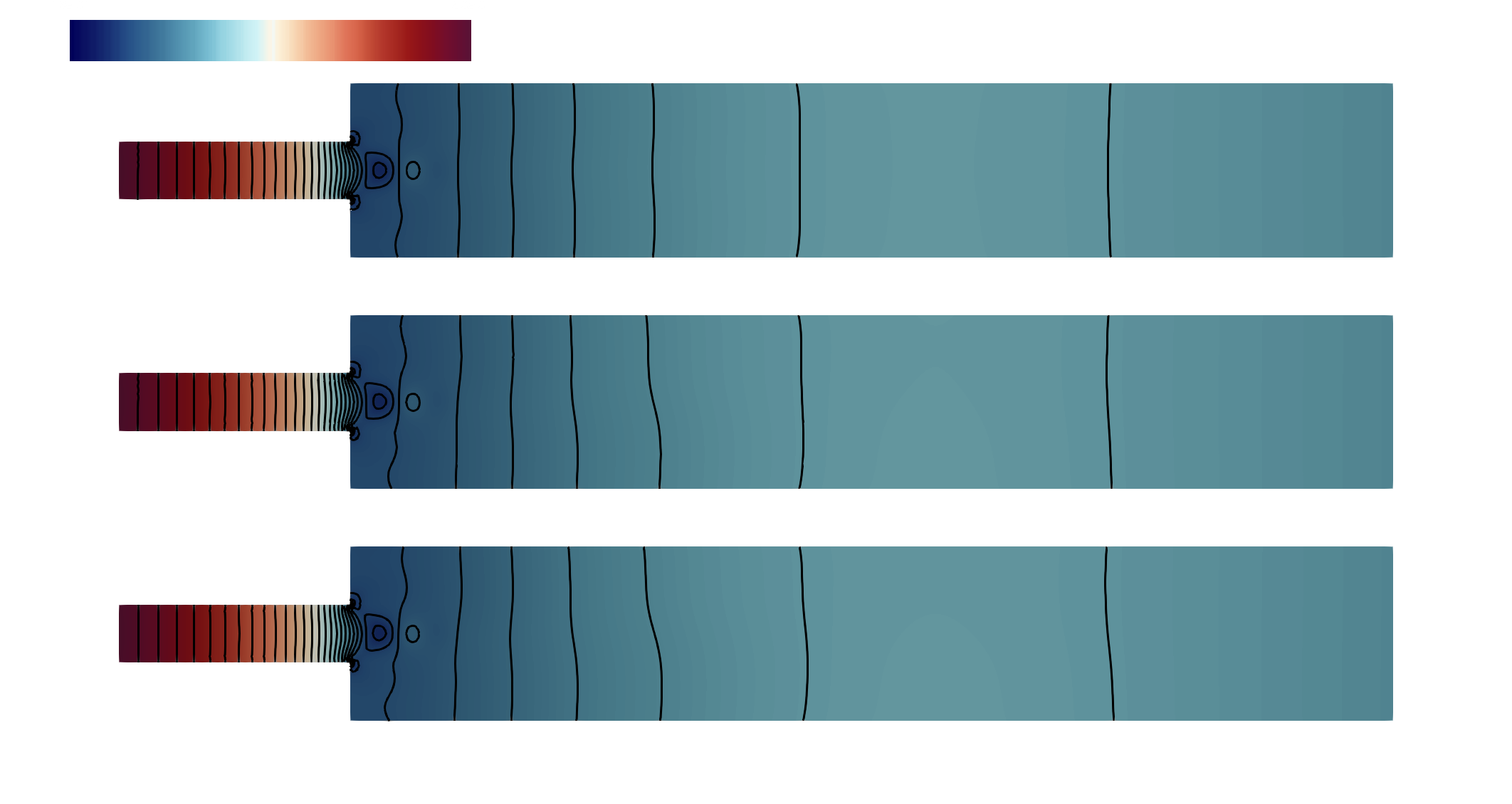}
\put (18.5,67.5) {$p/p_{\textrm{out}}$}
\put (2,65) {$0.13$}
\put (35,65) {$2.70$}
\put (0,50.7) {FOM}
\put (-9,31.8) {ROM RBF}
\put (-9,13.0) {ROM ANN}
\end{overpic}
\caption{Comparison of the pressure field between full-order model and different reduced-order models (for $\mu=0.77$  and $\textrm{Ma}=0.7$ ). Top, full order model; middle, reduced-order model using RBF interpolation; bottom, reduced-order model using ANN interpolation.}
\label{fig:comparison_mu077p}
\end{figure}
\begin{figure}
\centering
\begin{overpic}[trim=0 0 400 0, clip,width=0.9\textwidth]{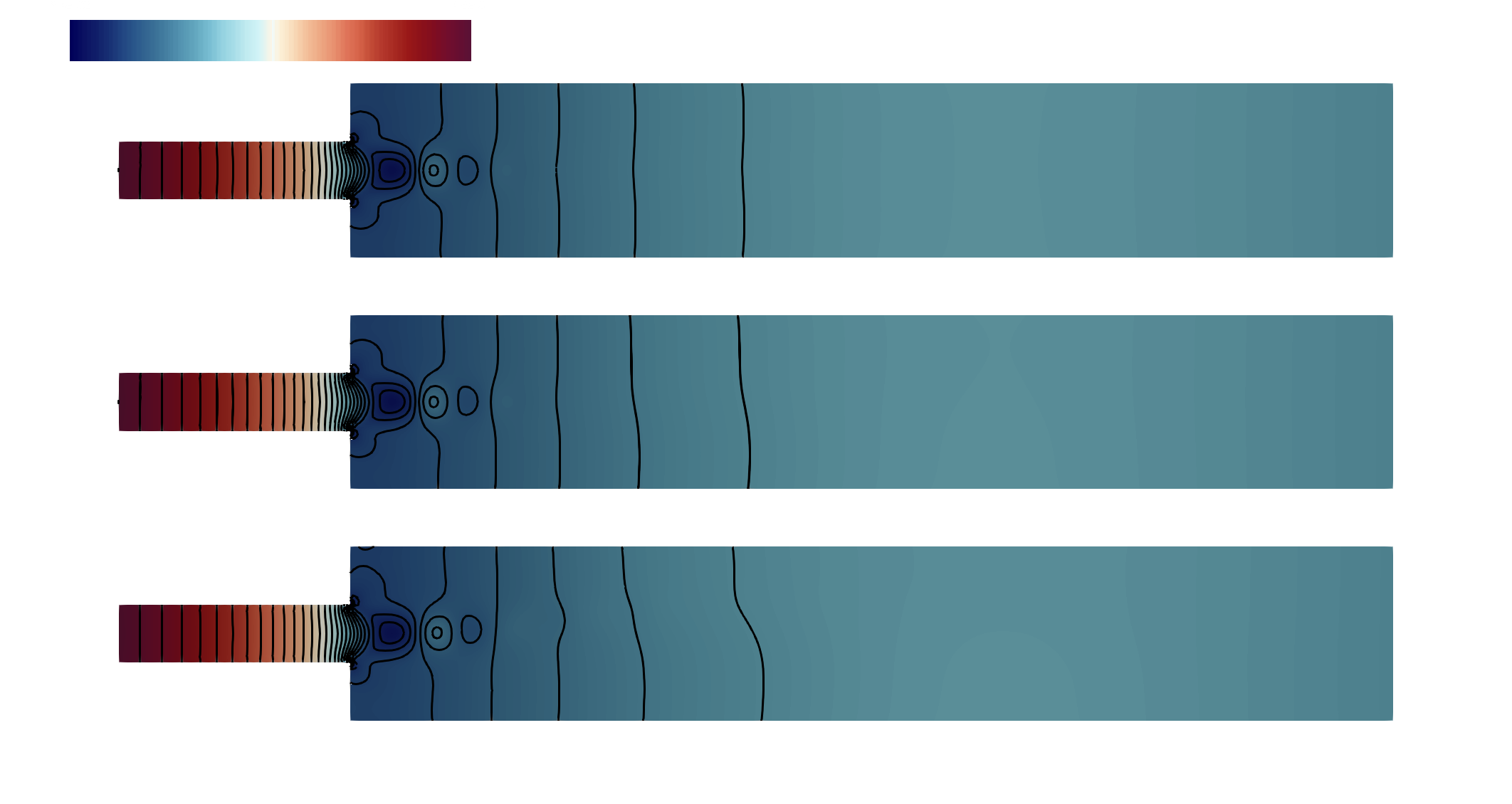}
\put (18.5,67.5) {$p/p_{\textbf{out}}$}
\put (2,65) {$0.094$}
\put (35,65) {$3.200$}
\put (0,50.7) {FOM}
\put (-9,31.8) {ROM RBF}
\put (-9,13.0) {ROM ANN}
\end{overpic}
\caption{Comparison of the pressure field between full-order model and different reduced-order models (for $\mu=0.67$  and $\textrm{Ma}=0.8$). Top, full order model; middle, reduced-order model using RBF interpolation; bottom, reduced-order model using ANN interpolation.}
\label{fig:comparison_mu067p}
\end{figure}
\begin{figure}
\centering
\begin{overpic}[trim=0 0 400 0, clip,width=0.9\textwidth]{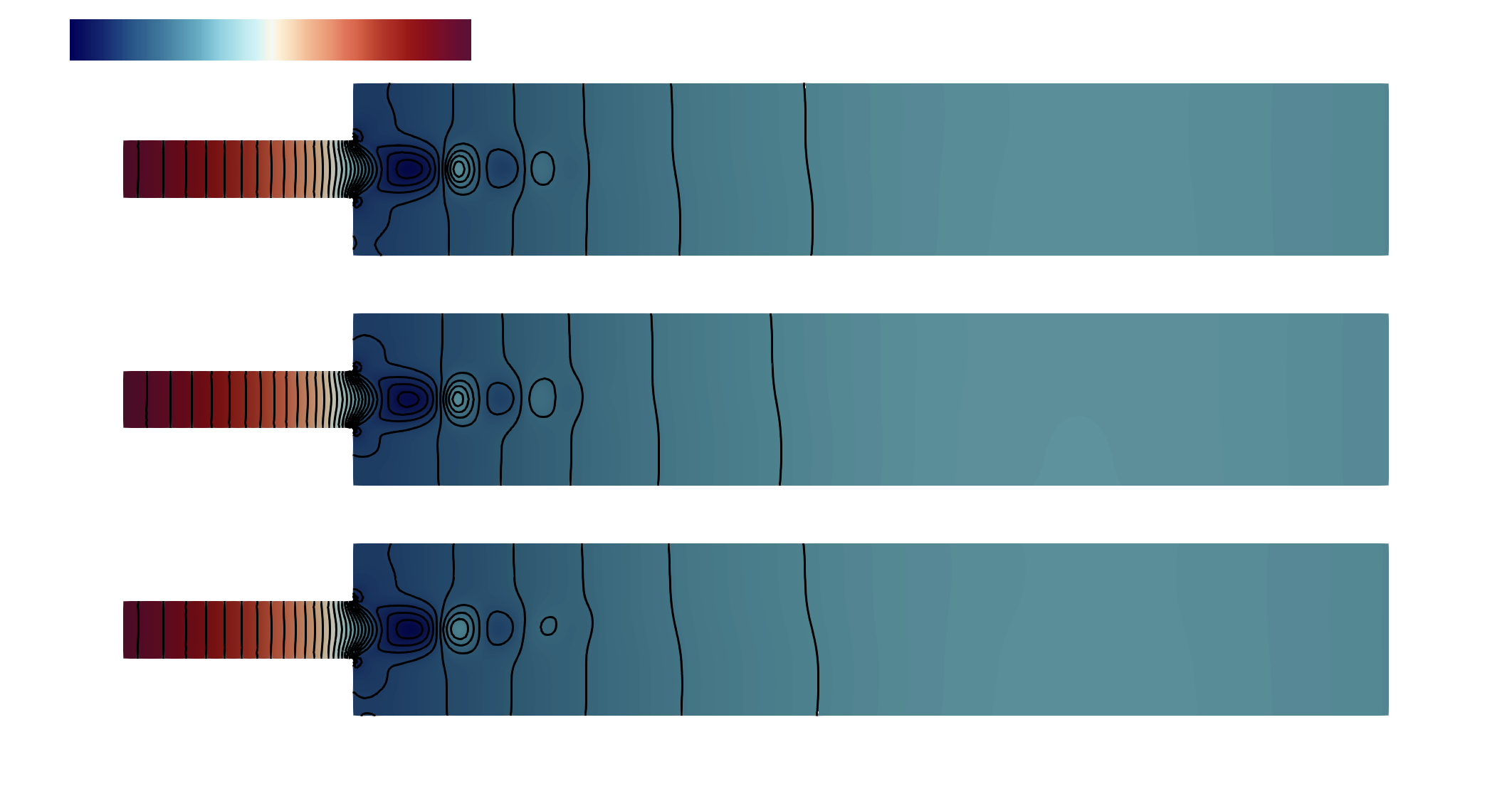}
\put (18.5,67.5) {$p/p_{\textbf{out}}$}
\put (2,65) {$0.062$}
\put (35,65) {$3.700$}
\put (0,50.7) {FOM}
\put (-9,31.8) {ROM RBF}
\put (-9,13.0) {ROM ANN}
\end{overpic}
\caption{Comparison of the pressure field between full-order model and different reduced-order models (for $\mu=0.55$  and $\textrm{Ma}=0.9$). Top, full order model; middle, reduced-order model using RBF interpolation; bottom, reduced-order model using ANN interpolation.}
\label{fig:comparison_mu055p}
\end{figure}

\section{Conclusions}\label{sec:5}
In the present work we employed two different non-intrusive Reduced Order Models (based on RBF and ANN interpolation) for the prediction of bifurcating phenomena in compressible fluid dynamics problems. In particular, a fully-compressible version of the classical Coanda effect in a suddenly expanding channel was considered as baseline configuration. A first part of the work was dedicated to a physical analysis of the influence of compressibility effects on the overall behaviour of the bifurcation. It was observed that for higher Mach numbers the stability region was extended allowing symmetric solutions for lower values of viscosity. In the second part of the paper, the two ROMs were employed to study the numerical modelling of the bifurcation. Both approaches were able to accurately predict the location of the bifurcation in terms of both velocity and pressure fields. Future work will be focused on highly compressible cases where shock waves naturally emerge. In these conditions, the development of more sophisticated reduced order models, based for example on deep neural network for learning non-linear regression, is mandatory in order to guarantee the consistency the reconstructed solution with its full order counterpart.

\section*{Acknowledgements}
This work was partially funded by European Union Funding for Research and Innovation --- Horizon 2020 Program --- in the framework of European Research
Council Executive Agency: H2020 ERC CoG 2015 AROMA-CFD project 681447 ``Advanced Reduced Order Methods with Applications in Computational Fluid Dynamics'', P.I. Professor Gianluigi Rozza.
The authors acknowledge support from MIUR (Italian ministry for university and research) through FARE-X-AROMA-CFD project, P.I. Prof. Gianluigi Rozza. Gianmarco Mengaldo acknowledges support from NUS startup grant R-265-000-A36-133.

\bibliography{./extracted.bib}

\begin{thebibliography}{10}
\expandafter\ifx\csname url\endcsname\relax
  \def\url#1{\texttt{#1}}\fi
\expandafter\ifx\csname urlprefix\endcsname\relax\def\urlprefix{URL }\fi
\expandafter\ifx\csname href\endcsname\relax
  \def\href#1#2{#2} \def\path#1{#1}\fi

\bibitem{slotnick2014cfd}
J.~P. Slotnick, A.~Khodadoust, J.~Alonso, D.~Darmofal, W.~Gropp, E.~Lurie,
  D.~J. Mavriplis, {CFD vision 2030 study: a path to revolutionary
  computational aerosciences}, Tech. rep. (2014).

\bibitem{hesthaven2007nodal}
J.~S. Hesthaven, T.~Warburton, {Nodal discontinuous Galerkin methods:
  algorithms, analysis, and applications}, Springer Science \& Business Media,
  2007.

\bibitem{cockburn:98}
B.~Cockburn, C.~Shu, The local discontinuous {G}alerkin finite element method
  for convection-diffusion systems, SIAM Journal of Numerical Analysis 35
  (1998) 2440--2463.

\bibitem{cockburn:98b}
B.~Cockburn, C.~Shu, The {R}unge-{K}utta discontinuous {G}alerkin finite
  element method for conservation laws {V}: Multidimensional systems, Journal
  of Computational Physics 141 (1998) 199--224.

\bibitem{shu1988efficient}
C.-W. Shu, S.~Osher, Efficient implementation of essentially non-oscillatory
  shock-capturing schemes, Journal of Computational Physics 77~(2) (1988)
  439--471.

\bibitem{liu1994weighted}
X.-D. Liu, S.~Osher, T.~Chan, Weighted essentially non-oscillatory schemes,
  Journal of Computational Physics 115~(1) (1994) 200--212.

\bibitem{shu1998essentially}
C.-W. Shu, Essentially non-oscillatory and weighted essentially non-oscillatory
  schemes for hyperbolic conservation laws, Advanced Numerical Approximation of
  Nonlinear Hyperbolic Equations (1998) 325--432.

\bibitem{Zhang2011PositivitypreservingHO}
X.~Zhang, C.-W. Shu, Positivity-preserving high order discontinuous {G}alerkin
  schemes for compressible {E}uler equations with source terms, Journal of
  Computational Physics 230 (2011) 1238--1248.

\bibitem{persson2013shock}
P.-O. Persson, {Shock Capturing for High-Order Discontinuous {G}alerkin
  Simulation of Transient Flow Problems}, AIAA Paper 2013-3061 (2013) 1--9,
  21st AIAA Computational Fluid Dynamics Conference, San Diego, CA,
  Jun.~24--27, 2013.

\bibitem{mengaldo2015discontinuous}
G.~Mengaldo, Discontinuous spectral/hp element methods: development, analysis
  and applications to compressible flows, Ph.D. thesis, Imperial College London
  (2015).

\bibitem{hillewaert2016assessment}
K.~Hillewaert, J.~Cagnone, S.~Murman, A.~Garai, Y.~Lv, M.~Ihme, Assessment of
  high-order {DG} methods for {LES} of compressible flows, in: Proceedings of
  the Summer Program, Center for Turbulence Research. Stanford University,
  2016, pp. 363--72.

\bibitem{moura2017setting}
R.~C. Moura, G.~Mengaldo, J.~Peir{\'o}, S.~J. Sherwin, An {LES} setting for
  {DG}-based implicit {LES} with insights on dissipation and robustness, in:
  Spectral and High Order Methods for Partial Differential Equations ICOSAHOM
  2016, Springer, Cham, 2017, pp. 161--173.

\bibitem{moxey2017towards}
D.~Moxey, C.~Cantwell, G.~Mengaldo, D.~Serson, D.~Ekelschot, J.~Peir{\'o},
  S.~Sherwin, R.~Kirby, Towards p-adaptive spectral/hp element methods for
  modelling industrial flows, in: Spectral and high order methods for partial
  differential equations icosahom 2016, Springer, Cham, 2017, pp. 63--79.

\bibitem{moura2017eddy}
R.~C. Moura, G.~Mengaldo, J.~Peir{\'o}, S.~J. Sherwin, On the eddy-resolving
  capability of high-order discontinuous {G}alerkin approaches to implicit
  {LES}/under-resolved {DNS} of {E}uler turbulence, Journal of Computational
  Physics 330 (2017) 615--623.

\bibitem{mengaldo2018spatial_1}
G.~Mengaldo, R.~Moura, B.~Giralda, J.~Peir{\'o}, S.~Sherwin, {Spatial
  eigensolution analysis of discontinuous Galerkin schemes with practical
  insights for under-resolved computations and implicit {LES}}, Computers \&
  Fluids 169 (2018) 349--364.

\bibitem{mengaldo2018spatial_2}
G.~Mengaldo, D.~De~Grazia, R.~C. Moura, S.~J. Sherwin, Spatial eigensolution
  analysis of energy-stable flux reconstruction schemes and influence of the
  numerical flux on accuracy and robustness, Journal of Computational Physics
  358 (2018) 1--20.

\bibitem{fernandez2018ability}
P.~Fernandez, N.-C. Nguyen, J.~Peraire, On the ability of discontinuous
  {G}alerkin methods to simulate under-resolved turbulent flows, arXiv preprint
  arXiv:1810.09435 (2018).

\bibitem{winters2018comparative}
A.~R. Winters, R.~C. Moura, G.~Mengaldo, G.~J. Gassner, S.~Walch, J.~Peiro,
  S.~J. Sherwin, A comparative study on polynomial dealiasing and split form
  discontinuous {G}alerkin schemes for under-resolved turbulence computations,
  Journal of Computational Physics 372 (2018) 1--21.

\bibitem{fernandez2019non}
P.~Fernandez, R.~C. Moura, G.~Mengaldo, J.~Peraire, Non-modal analysis of
  spectral element methods: {T}owards accurate and robust large-eddy
  simulations, Computer Methods in Applied Mechanics and Engineering 346 (2019)
  43--62.

\bibitem{moura2020viscous}
R.~C. Moura, P.~Fernandez, G.~Mengaldo, S.~J. Sherwin, Viscous diffusion
  effects in the eigenanalysis of (hybridisable) {DG} methods, in: Spectral and
  High Order Methods for Partial Differential Equations ICOSAHOM 2018,
  Springer, Cham, 2020, pp. 371--382.

\bibitem{mengaldo2021industry}
G.~Mengaldo, D.~Moxey, M.~Turner, R.~C. Moura, A.~Jassim, M.~Taylor, J.~Peiro,
  S.~Sherwin, Industry-relevant implicit large-eddy simulation of a
  high-performance road car via spectral/hp element methods, SIAM Review 63~(4)
  (2021) 723--755.

\bibitem{moura2022spectral}
R.~Moura, L.~Fernandes, A.~Silva, G.~Mengaldo, S.~Sherwin, Spectral/hp element
  methods' linear mechanism of (apparent) energy transfer in {F}ourier space:
  Insights into dispersion analysis for implicit {LES}, Journal of
  Computational Physics 471 (2022) 111613.

\bibitem{hesthaven2016certified}
J.~S. Hesthaven, G.~Rozza, B.~Stamm, et~al., Certified reduced basis methods
  for parametrized partial differential equations, Vol. 590, Springer, 2016.

\bibitem{quarteroni2015reduced}
A.~Quarteroni, A.~Manzoni, F.~Negri, Reduced basis methods for partial
  differential equations: an introduction, Vol.~92, Springer, 2015.

\bibitem{chinesta2016model}
F.~Chinesta, A.~Huerta, G.~Rozza, K.~Willcox, Model order reduction,
  Encyclopedia of computational mechanics (2016).

\bibitem{benner2017model}
P.~Benner, M.~Ohlberger, A.~Patera, G.~Rozza, K.~Urban, Model reduction of
  parametrized systems, Springer, 2017.

\bibitem{hotelling1933analysis}
H.~Hotelling, Analysis of a complex of statistical variables into principal
  components., Journal of Educational Psychology 24~(6) (1933) 417.

\bibitem{lumley1967structure}
J.~L. Lumley, The structure of inhomogeneous turbulent flows, Atmospheric
  Turbulence and Radio Wave Propagation (1967) 166--178.

\bibitem{burkardt2006pod}
J.~Burkardt, M.~Gunzburger, H.-C. Lee, Pod and cvt-based reduced-order modeling
  of navier--stokes flows, Computer Methods in Applied Mechanics and
  Engineering 196~(1-3) (2006) 337--355.

\bibitem{ballarin2015supremizer}
F.~Ballarin, A.~Manzoni, A.~Quarteroni, G.~Rozza, {Supremizer stabilization of
  POD--Galerkin approximation of parametrized steady incompressible
  Navier--Stokes equations}, International Journal for Numerical Methods in
  Engineering 102~(5) (2015) 1136--1161.

\bibitem{maulik2021pyparsvd}
R.~Maulik, G.~Mengaldo, {PyParSVD}: A streaming, distributed and randomized
  singular-value-decomposition library, in: 2021 7th International Workshop on
  Data Analysis and Reduction for Big Scientific Data (DRBSD-7), IEEE, 2021,
  pp. 19--25.

\bibitem{hess2022reduced}
M.~W. Hess, A.~Lario, G.~Mengaldo, G.~Rozza, Reduced order modeling for
  spectral element methods: current developments in {Nektar++} and further
  perspectives, arXiv preprint arXiv:2201.05404 (2022).

\bibitem{towne2018spectral}
A.~Towne, O.~T. Schmidt, T.~Colonius, Spectral proper orthogonal decomposition
  and its relationship to dynamic mode decomposition and resolvent analysis,
  Journal of Fluid Mechanics 847 (2018) 821--867.

\bibitem{schmidt2019spectral}
O.~T. Schmidt, G.~Mengaldo, G.~Balsamo, N.~P. Wedi, Spectral empirical
  orthogonal function analysis of weather and climate data, Monthly Weather
  Review 147~(8) (2019) 2979--2995.

\bibitem{mengaldo2021pyspod}
G.~Mengaldo, R.~Maulik, {PySPOD}: A python package for spectral proper
  orthogonal decomposition ({SPOD}), Journal of Open Source Software 6~(60)
  (2021) 2862.

\bibitem{lario2022neural}
A.~Lario, R.~Maulik, O.~T. Schmidt, G.~Rozza, G.~Mengaldo, Neural-network
  learning of {SPOD} latent dynamics, Journal of Computational Physics 468
  (2022) 111475.

\bibitem{yano2019}
M.~Yano, {Discontinuous Galerkin reduced basis empirical quadrature procedure
  for model reduction of parametrized nonlinear conservation laws}, Advances in
  Computational Mathematics 45 (12 2019).

\bibitem{BALAJEWICZ2016224}
M.~Balajewicz, I.~Tezaur, E.~Dowell, {Minimal subspace rotation on the Stiefel
  manifold for stabilization and enhancement of projection-based reduced order
  models for the compressible Navier–Stokes equations}, Journal of
  Computational Physics 321 (2016) 224--241.

\bibitem{WU2020112766}
P.~Wu, J.~Sun, X.~Chang, W.~Zhang, R.~Arcucci, Y.~Guo, C.~C. Pain, Data-driven
  reduced order model with temporal convolutional neural network, Computer
  Methods in Applied Mechanics and Engineering 360 (2020) 112766.

\bibitem{pichi2020reduced}
F.~Pichi, J.~Eftang, G.~Rozza, A.~Patera, Reduced order models for the buckling
  of hyperelastic beams (2020).

\bibitem{pichi2019reduced}
F.~Pichi, G.~Rozza, {Reduced basis approaches for parametrized bifurcation
  problems held by non-linear Von K{\'a}rm{\'a}n equations}, Journal of
  Scientific Computing 81~(1) (2019) 112--135.

\bibitem{pichi2022driving}
F.~Pichi, M.~Strazzullo, F.~Ballarin, G.~Rozza, {Driving bifurcating
  parametrized nonlinear PDEs by optimal control strategies: application to
  Navier--Stokes equations with model order reduction}, ESAIM: Mathematical
  Modelling and Numerical Analysis 56~(4) (2022) 1361--1400.

\bibitem{pitton2017computational}
G.~Pitton, A.~Quaini, G.~Rozza, Computational reduction strategies for the
  detection of steady bifurcations in incompressible fluid-dynamics:
  Applications to coanda effect in cardiology, Journal of Computational Physics
  344 (2017) 534--557.

\bibitem{pitton2017application}
G.~Pitton, G.~Rozza, On the application of reduced basis methods to bifurcation
  problems in incompressible fluid dynamics, Journal of Scientific Computing
  73~(1) (2017) 157--177.

\bibitem{martin2020reduced}
M.~W. Hess, A.~Quaini, G.~Rozza, {Reduced basis model order reduction for
  Navier–Stokes equations in domains with walls of varying curvature},
  International Journal of Computational Fluid Dynamics 34~(2) (2020) 119--126.

\bibitem{hess2022sparse}
M.~W. Hess, G.~Rozza, {Model Reduction Using Sparse Polynomial Interpolation
  for the Incompressible Navier–Stokes Equations}, Vietnam Journal of
  Mathematics (2022).

\bibitem{khamlich2021model}
M.~Khamlich, F.~Pichi, G.~Rozza, {Model order reduction for bifurcating
  phenomena in Fluid-Structure Interaction problems}, arXiv preprint
  arXiv:2110.06297 (2021).

\bibitem{tritton2012physical}
D.~J. Tritton, Physical fluid dynamics, Springer Science \& Business Media,
  2012.

\bibitem{ahmed2019coanda}
N.~A. Ahmed, Coanda Effect: flow phenomenon and applications, CRC Press, 2019.

\bibitem{drikakis1997bifurcation}
D.~Drikakis, Bifurcation phenomena in incompressible sudden expansion flows,
  Physics of Fluids 9~(1) (1997) 76--87.

\bibitem{allery2004application}
C.~Allery, J.-M. Cadou, A.~Hamdouni, D.~Razafindralandy, {Application of the
  asymptotic numerical method to the Coanda effect study}, Revue Europ{\'e}enne
  des {\'E}l{\'e}ments 13~(1-2) (2004) 57--77.

\bibitem{saha2020bifurcation}
S.~Saha, P.~Biswas, S.~Nath, {Bifurcation phenomena for incompressible laminar
  flow in expansion channel to study Coanda effect}, Journal of
  interdisciplinary Mathematics 23~(2) (2020) 493--502.

\bibitem{haffner2020unsteady}
Y.~Haffner, J.~Bor{\'e}e, A.~Spohn, T.~Castelain, {Unsteady Coanda effect and
  drag reduction for a turbulent wake}, Journal of Fluid Mechanics 899 (2020).

\bibitem{cantwell2015nektar++}
C.~D. Cantwell, D.~Moxey, A.~Comerford, A.~Bolis, G.~Rocco, G.~Mengaldo,
  D.~De~Grazia, S.~Yakovlev, J.-E. Lombard, D.~Ekelschot, et~al., {Nektar++: An
  open-source spectral/hp element framework}, Computer Physics Communications
  192 (2015) 205--219.

\bibitem{moxey2020nektar++}
D.~Moxey, C.~D. Cantwell, Y.~Bao, A.~Cassinelli, G.~Castiglioni, S.~Chun,
  E.~Juda, E.~Kazemi, K.~Lackhove, J.~Marcon, et~al., {Nektar++: Enhancing the
  capability and application of high-fidelity spectral/hp element methods},
  Computer Physics Communications 249 (2020) 107110.

\bibitem{mengaldo2014guide}
G.~Mengaldo, D.~De~Grazia, J.~Peiro, A.~Farrington, F.~Witherden, P.~Vincent,
  S.~Sherwin, {A Guide to the Implementation of Boundary Conditions in Compact
  High-Order Methods for Compressible Aerodynamics}, in: AIAA Aviation Forum
  2014, 2014.

\bibitem{mengaldo2015dealiasing}
G.~Mengaldo, D.~De~Grazia, D.~Moxey, P.~E. Vincent, S.~J. Sherwin, Dealiasing
  techniques for high-order spectral element methods on regular and irregular
  grids, Journal of Computational Physics 299 (2015) 56--81.

\bibitem{mengaldo2015triple}
G.~Mengaldo, M.~Kravtsova, A.~Ruban, S.~Sherwin, Triple-deck and direct
  numerical simulation analyses of high-speed subsonic flows past a roughness
  element, Journal of Fluid Mechanics 774 (2015) 311--323.

\bibitem{lombard2016implicit}
J.-E.~W. Lombard, D.~Moxey, S.~J. Sherwin, J.~F. Hoessler, S.~Dhandapani, M.~J.
  Taylor, Implicit large-eddy simulation of a wingtip vortex, AIAA Journal
  54~(2) (2016) 506--518.

\bibitem{serson2017direct}
D.~Serson, J.~R. Meneghini, S.~J. Sherwin, Direct numerical simulations of the
  flow around wings with spanwise waviness, Journal of Fluid Mechanics 826
  (2017) 714--731.

\bibitem{nakhchi2020dns}
M.~E. Nakhchi, S.~W. Naung, M.~Rahmati, {DNS} of secondary flows over
  oscillating low-pressure turbine using spectral/hp element method,
  International Journal of Heat and Fluid Flow 86 (2020) 108684.

\bibitem{barone2009stable}
M.~F. Barone, I.~Kalashnikova, D.~J. Segalman, H.~K. Thornquist, {Stable
  Galerkin reduced order models for linearized compressible flow}, Journal of
  Computational Physics 228~(6) (2009) 1932--1946.

\bibitem{mcquarrie2021data}
S.~A. McQuarrie, C.~Huang, K.~E. Willcox, Data-driven reduced-order models via
  regularised operator inference for a single-injector combustion process,
  Journal of the Royal Society of New Zealand 51~(2) (2021) 194--211.

\bibitem{wang2009transient}
T.-S. Wang, Transient three-dimensional startup side load analysis of a
  regeneratively cooled nozzle, Shock Waves 19~(3) (2009) 251--264.

\bibitem{trancossi2011overview}
M.~Trancossi, {An overview of scientific and technical literature on Coanda
  effect applied to nozzles} (2011).

\bibitem{trancossi2011acheon}
M.~Trancossi, A.~Dumas, {ACHEON: Aerial Coanda High Efficiency Orienting-jet
  Nozzle}, Tech. rep., SAE Technical Paper (2011).

\bibitem{ahmed2017aerodynamics}
R.~Ahmed, A.~A. Talib, A.~M. Rafie, H.~Djojodihardjo, Aerodynamics and flight
  mechanics of mav based on coanda effect, Aerospace Science and Technology 62
  (2017) 136--147.

\bibitem{lubert2010some}
C.~Lubert, {On some recent applications of the Coanda effect to acoustics}, in:
  Proceedings of Meetings on Acoustics 160ASA, Vol.~11, Acoustical Society of
  America, 2010, p. 040006.

\bibitem{freire2002bubble}
A.~P.~S. Freire, D.~D. Miranda, L.~M. Luz, G.~F. Fran{\c{c}}a, {Bubble plumes
  and the Coanda effect}, International Journal of Multiphase Flow 28~(8)
  (2002) 1293--1310.

\bibitem{ginghina2007coandua}
C.~Ginghina, {The Coanda effect in cardiology}, Journal of Cardiovascular
  Medicine 8~(6) (2007) 411--413.

\bibitem{hess2019localized}
M.~Hess, A.~Alla, A.~Quaini, G.~Rozza, M.~Gunzburger, {A localized
  reduced-order modeling approach for PDEs with bifurcating solutions},
  Computer Methods in Applied Mechanics and Engineering 351 (2019) 379--403.

\bibitem{pintore2021efficient}
M.~Pintore, F.~Pichi, M.~Hess, G.~Rozza, C.~Canuto, Efficient computation of
  bifurcation diagrams with a deflated approach to reduced basis spectral
  element method, Advances in Computational Mathematics 47~(1) (2021) 1--39.

\bibitem{pichi2021artificial}
F.~Pichi, F.~Ballarin, G.~Rozza, J.~S. Hesthaven, An artificial neural network
  approach to bifurcating phenomena in computational fluid dynamics, arXiv
  preprint arXiv:2109.10765 (2021).

\bibitem{durst1974low}
F.~Durst, A.~Melling, J.~H. Whitelaw, {Low Reynolds number flow over a plane
  symmetric sudden expansion}, Journal of Fluid Mechanics 64~(1) (1974)
  111--128.

\bibitem{fearn1990nonlinear}
R.~Fearn, T.~Mullin, K.~Cliffe, Nonlinear flow phenomena in a symmetric sudden
  expansion, Journal of Fluid Mechanics 211 (1990) 595--608.

\bibitem{karantonis2021compressibility}
K.~Karantonis, I.~W. Kokkinakis, B.~Thornber, D.~Drikakis, Compressibility in
  suddenly expanded subsonic flows, Physics of Fluids 33~(10) (2021) 105106.

\bibitem{kennedy2016diagonally}
C.~A. Kennedy, M.~H. Carpenter, Diagonally implicit runge-kutta methods for
  ordinary differential equations. a review, Tech. rep. (2016).

\bibitem{toro2013riemann}
E.~F. Toro, Riemann solvers and numerical methods for fluid dynamics: a
  practical introduction, Springer Science \& Business Media, 2013.

\bibitem{hartmann2006symmetric}
R.~Hartmann, P.~Houston, Symmetric interior penalty dg methods for the
  compressible navier--stokes equations ii: Goal--oriented a posteriori error
  estimation, International Journal of Numerical Analysis \& Modeling 3~(2)
  (2006) 141--162.

\bibitem{cagniart2019model}
N.~Cagniart, Y.~Maday, B.~Stamm, Model order reduction for problems with large
  convection effects, in: Contributions to partial differential equations and
  applications, Springer, 2019, pp. 131--150.

\bibitem{demo18ezyrb}
N.~Demo, M.~Tezzele, G.~Rozza, {EZyRB: Easy Reduced Basis method}, The Journal
  of Open Source Software 3~(24) (2018) 661.

\bibitem{nair2019transported}
N.~J. Nair, M.~Balajewicz, Transported snapshot model order reduction approach
  for parametric, steady-state fluid flows containing parameter-dependent
  shocks, International Journal for Numerical Methods in Engineering 117~(12)
  (2019) 1234--1262.

\end{thebibliography}
\end{document}